\documentclass[12pt]{article}
\usepackage{amssymb}
\usepackage{latexsym}
\usepackage[T1]{fontenc}
\usepackage[cp1250]{inputenc}
\newcounter{myownsection}
\def\myownsection{\refstepcounter{myownsection} \setcounter{equation}{0}}
\usepackage[dvips]{graphicx}
\usepackage{amsmath}
\usepackage{amssymb}
\usepackage{amsfonts}
\usepackage{eufrak}
\usepackage{fancyhdr}
\usepackage[mathscr]{eucal}

\begin{document}
\begin{center}
{\bf DISCRETE INTERPOLATION BETWEEN MONOTONE\\
PROBABILITY AND FREE PROBABILITY}\footnote{This work is supported by KBN grant No 2P03A00723 and by the EU Network QP-Applications,
Contract No. HPRN-CT-2002-00279}\\[60pt]
{\bf Romuald Lenczewski and Rafa\l{} Sa\l{}apata}\\[40pt]
{\it Institute of Mathematics\\
Wroc\l{}aw University of Technology\\
Wybrze\.{z}e Wyspia\'nskiego 27\\
50-370 Wroc\l{}aw\\
Poland}\\[20pt]
\end{center}
\begin{abstract}
We construct a sequence of states called {\it $m$-monotone product states} which give
a discrete interpolation between the monotone product of states of Muraki and the free product of states
of Avitzour and Voiculescu in free probability. We derive the associated basic limit theorems
and develop the combinatorics based on non-crossing ordered partitions with
monotone order starting from depth $m$. The Hilbert space representations of the limit mixed moments in
the invariance principle lead to $m$-{\it monotone Gaussian operators} living in $m$-{\it monotone Fock spaces},
which are truncations of the free Fock space ${\cal F}(L^{2}({\mathbb R}^{+}))$.
A new type of combinatorics of inner blocks in non-crossing partitions leads
to explicit formulas for the mixed moments of $m$-monotone Gaussian operators,
which are new even in the case of monotone independent Gaussian operators with arcsine
distributions.  \\[10pt]
\end{abstract}

\myownsection
\begin{center}
{\sc 1. Introduction.}
\end{center}
The aim of this paper is to construct an interpolation between free probability [V,VDN] based
on the notion of the free product of states [Av,V] and monotone probability
based on the notion of the monotone product of states [M2]. Our interpolation is given by
a sequence of product states called $m$-{\it monotone}, where $m\in {\mathbb N}$, which
gives the monotone product for $m=1$, whose asymptotics leads to the free product for $m=\infty$.

The definition of the free product of states is based on the freeness condition
\begin{equation}
\phi(a_{1}a_{2}...a_{n})=0\;\;\; {\rm for}\;\;
a_{k}\in{\mathcal{A}}_{i_k} \cap {\rm Ker}{\phi}_{i_k}
\end{equation}
on the free product $*_{i\in I}{\cal A}_{i}$ of unital *-algebras with identified units
($1_{i}=1_{k}$ for all $i,k\in I$), where $i_1\neq i_2\neq \ldots \neq i_n$.

It is interesting to observe that the monotone product of states can also be defined
using (1.1), except that one has to take the free product without identification of units
$\sqcup_{i\in I}{\cal A}_{i}$, provided $I$ is linearly ordered, and give conditions on moments involving
units, namely
\begin{equation}
\phi(a_1 ... a_{k-1}1_{i_k}a_{k+1}...a_n)
=
\phi(a_1 ... a_{k-1}a_{k+1}...a_n)
\end{equation}
for $i_1<i_2<\ldots <i_k$ and $a_{1}, \ldots, a_{k-1}$ taken from the kernels of the states
(if the first unit is in a different configuration of indices, then it vanishes).

This observation leads to a construction of a sequence of states
$(\phi^{(m)})_{m=1}^{\infty}$ on $\sqcup_{i\in I}{\cal A}_{i}$, defined
by (1.1) on $\sqcup_{i\in I}{\cal A}_{i}$ and identifying
units on an increasing sequence of `subspaces'. This gives a discrete interpolation between the
monotone product for $m=1$ and the free product obtained in the limit $m\rightarrow \infty$,
which in some sense `neglects order' with increasing $m$.
The sequence of product states is called the {\it monotone hierarchy}
and the associated `independent' variables are called $m$-{\it monotone}.

The pace of this interpolation is based on the combinatorics of ordered non-crossing partitions
with monotone order (color) starting from depth $m$, giving monotone non-crossing partitions for
$m=1$ and all ordered non-crossing partitions for $m=\infty$.
In particular, the (even) moments of the standard (mean zero, variance one) central limit laws are
given by
\begin{equation}
M_{2k}^{(m)}
=
1/k!|{\cal ONC}_{2k}^{2}(m)|
\end{equation}
where ${\cal ONC}_{2k}^{2}(m)$ stands for ordered non-crossing pair partitions of a
$2k$-element set, which
have a monotone order starting from depth $m$.
In particular, for $m=1$ we get the moments $M_{2k}^{(1)}
=2^{-k}{2k \choose k}$ of the arcsine law
supported on $[-\sqrt{2}, \sqrt{2}]$,
and in the limit $m\rightarrow \infty$ we obtain
the moments $M_{2k}^{(\infty)}=\frac{1}{k+1}{2k \choose k}$ of the
Wigner measure supported on $[-2,2]$.

On the level of Cauchy transforms of central limit laws, the interpolation can be
nicely described in terms of their continued fractions.
Thus, let
\begin{equation}
G(z)= \cfrac{1}{z-\cfrac{\beta_{1}}{z-\cfrac{\beta_{2}}{z-\cfrac{\beta_{3}}{\ddots}}}}\
\end{equation}
be the Cauchy transform of a probability measure associated
with the Jacobi sequence $(\beta_{n})=(\beta_{1},\beta_{2},\beta_{3}, \ldots )$.
Then, to sequences
\begin{equation}
(\beta_{n}^{(m)})=(\underbrace{1,1,\ldots , 1}_{m \;{\rm times}},1/2,1/2, \ldots )
\end{equation}
correspond Cauchy transforms $G^{(m)}(z)$ of the $m$-monotone central limit measures.
In particular, for $m=1$ we get
$(\beta_{n}^{(1)})=(1,1/2,1/2, \ldots )$, corresponding to the arcsine law,
and for $m=\infty$, we get $(\beta_{n}^{(\infty)})=(1,1,1, \ldots)$, corresponding to
the Wigner law.
Thus, succesive iterations of $G^{(\infty)}(z)$ are obtained from
$G^{(1)}(z)$ by taking the $m-1$-th power of the right-sided shift
on the set of Jacobi sequences of the form (1.3), with $1$'s replacing $1/2$'s at the first $m$ sites.

We also derive the invariance principle and give a realization of
the limits using $m$-{\it monotone Gaussian operators} $\omega^{(m)}(f)$, where
$f\in \Theta=\{\chi_{(s,t]}: 0\leq s<t\}$, on truncations ${\cal F}^{(m)}({\cal H})$
of the free Fock space ${\cal F}({\cal H})$, where ${\cal H}=L^{2}({\mathbb R}^{+})$,
called $m$-{\it monotone Fock spaces}.
If $f_{1}, f_{2}, \ldots , f_{2k}\in \Theta$ have identical or disjoint supports, then
explicit combinatorial formulas for the mixed moments of
$\omega^{(m)}(f_{1})$, $\omega^{(m)}(f_{2}), \ldots, \omega^{(m)}(f_{2k})$ in the vacuum state
$\varphi$ are derived.
It is worth noting that even in the monotone case we get a new formula
\begin{equation}
\varphi(\omega^{(1)}(f_{1})\omega^{(1)}(f_{2})\ldots \omega^{(1)}(f_{2k}))
=
\sum_{
\stackrel{\pi\in {\cal NC}_{2k}^{2}}
{\scriptscriptstyle \pi \sim (f_1, \ldots ,f_{2k})}
}
\prod_{i=1}^{k}\frac{t^{(i)}-s^{(i)}}{{\rm Inn}(\pi_{i})+1}
\end{equation}
where summation runs over non-crossing partitions $\pi=\{\pi_{1}, \pi_{2}, \ldots , \pi_{k}\}$
`compatible' (in the sense specified in Section 10) with the supports of $f_{1},f_{2}, \ldots , f_{n}$
(support $(s^{(i)},t^{(i)}]$ is associated with block $\pi_{i}$),
with ${\rm Inn}(\pi_{i})$ denoting the number of blocks which are inner with respect to
block $\pi_{i}$ and are associated with supports identical to that associated with $\pi_{i}$.

Our approach resembles the construction of the {\it hierarchy of freeness} [L1] (see [FL] for the limit
theorems) which gave a discrete interpolation between the boolean product of states and the
free product of states. The latter is understood in [L1] in the weak sense of convergence of moments.
However, a reduction of free independence (free product of states) to tensor independence
(tensor product of states) in the strong sense (reduction of the free product
of states to the tensor product of states on extended algebras) has recently also been given [L2].
A similar approach to that in [L2] can also be given for the monotone hierarchy.

Other interpolations between various models of quantum probability
were studied in [BLS], [BW]. \\[10pt]

\myownsection
\begin{center}
{\sc 2. Products of states}
\end{center}
In this section we introduce a new definition of the monotone product of states, which
allows us to compare it with the free product of states and later leads to the construction
of a discrete interpolation between the two products.

By a {\it noncommutative probability space} we understand a pair \textup{($\mathcal{A}\,,\phi$)},
where $\mathcal{A}$ is a unital *-algebra over $\mathbb{C}$ and
$\phi\,:\mathcal{A}$ $\to$  $\mathbb{C}$ is a positive ($\phi(x^*x)\geqslant 0$ for all
$x\in\mathcal{A}$) normalized ($\phi(1)=1$) linear functional called {\it state}.
Elements of the *-algebra $\mathcal{A}$ are called {\it random variables}.
A noncommutative probability space \textup{($\mathcal{A}\,,\phi$)} is called
a $\mathcal{C}$*-{\it probability space} if $\mathcal{A}$ is a unital $\mathcal{C}$*-algebra.
By the {\it distribution} of the random variable $a \in\mathcal{A}$ we understand the linear functional
$\mu_{a}:\mathbb{C}[X]\to\mathbb{C}$, defined on the polynomials of $X$ by
\begin{center}
$\mu_{a}(\textup{P})=\phi\big(\textup{P}(a)\big)$.
\end{center}
Of course, the distribution $\mu_{a}$ is uniquely determined by the moments
$\phi(a^{n})$ of the variable $a$ for $n\in\mathbb{N}$.

We assume in this paper that $I$ is a linearly ordered set of indices.
\\
\indent{\par}
{\sc Definition 2.1.}
By the {\it monotone product} of states $\phi_{i}$ we understand the state
$\phi$ on the free product $\sqcup_{i\in I}{\cal A}_{i}$ without identification of units
given by
\begin{equation}
\phi(a_1 ... a_{k-1}a_{k}a_{k+1}...a_n)=\phi_{i_k}(a_k)\phi(a_1 ... a_{k-1}a_{k+1}...a_n),
\end{equation}
for $a_j\in\mathcal{A}_{i_j}$ whenever neighboring variables come from different algebras and
$i_{k-1} < i_k > i_{k+1}$ (with only one inequality if $k\in\{1,n\}$).\\
\indent{\par}
Note that the freeness condition (1.1) is of quite different form than
the monotonicity condition given by (2.1). However, as we demonstrate below,
one can reduce Definition 2.1 to the freeness condition plus
conditions on moments involving units.
In the case of the free product of states, the condition on moments involving units
is trivial, namely units of *-algebras ${\cal A}_{i}$ always play the role of a `global unit' and can be identified.
In the case of the monotone product of states, units of *-algebras ${\cal A}_{i}$ play the role of a
`monotone approximate unit'.
This can be done even if *-algebras ${\cal A}_{i}$ are not unital - it is enough to take their
unitizations. \\
\indent{\par}
{\sc Definition 2.2.}
Let $({\cal A}_{i},\phi_{i})_{i\in I}$ be a family of noncommutative probability spaces.
Define $\psi: \sqcup_{i\in I}{\cal A}_{i}\rightarrow {\mathbb C}$ to be the linear functional
defined by the freeness condition (1.1) for all $a_{j}$ in the kernels of the $\phi_{i_{j}}$, respectively,
and by the condition
\begin{equation}
\psi(a_1 ... a_{k-1}1_{i_k}a_{k+1}...a_n)
=
\left\{
\begin{array}{cc}
\psi(a_1 ... a_{k-1}a_{k+1}...a_n) & {\rm if}\;\;i_1<i_2<\ldots <i_k\\
0                                  & {\rm otherwise}
\end{array}
\right.
\end{equation}
where $1\leq k\leq n$ and only the variables $a_{j}$ preceeding $1_{i_{k}}$ are assumed to be in the
kernels of the states $\phi_{i_j}$, respectively, with $\phi_{i}=\psi|_{{\cal A}_{i}}$ for every $i\in I$.\\
\indent{\par}
{\it Example 2.1.}
Let $I=\{1,2\}$ (with natural order) and $a \in {\cal A}_{1}$, $b\in {\cal A}_{2}$ and denote
$a^{0}=a-\psi(a)1_{1}$, $b^{0}=b-\psi(b)1_{2}$.
Then $\psi(a^{0}1_2a^{0})=\psi((a^{0})^2)$ and
$\psi(b^{0}1_1b^{0})=0$
are the key formulas which, together with (1.1), lead to different factorizations
\begin{eqnarray*}
\psi(aba)&=&\psi(a^{2})\psi(b)\\
\psi(bab)&=&\psi^{2}(b)\psi(a)
\end{eqnarray*}
characteristic of the monotone product of states.\\
\indent{\par}
{\sc Theorem 2.3.}
{\it The state $\psi$ agrees with the monotone product of states $\phi_{i}$, $i\in I$.}\\
\indent{\par}
{\it Proof.}
Let $a_1\in\mathcal{A}_{i_1}$, $a_2\in\mathcal{A}_{i_2}$,\ldots , $a_n\in\mathcal{A}_{i_n}$.
Denote $a_j^0=a_j-\phi_{i_j}(a_j)1_{i_j}$ for any $a_j \in {\cal A}_{i_j}$ and let $\phi$
stand for the monotone product of states $\phi_{i}$.
Using (2.2) for $n=1$ and the fact that $\psi(1_{i_1})=1$, we obtain
$$
\psi(a_{1})=\psi(a_{1}^{0}+\phi_{i_{1}}(a_{1})1_{i_{1}})=\phi_{i_{1}}(a_{1}).
$$
Suppose now that $\psi$ agrees with the monotone product of $\phi_{i}$'s
on words of lenght smaller than $n$ and assume that
$i_1 < i_2 <...< i_k > i_{k+1}$ (if $k=n$, then the last inequality is ignored). Decomposing
$a_j=a_j^0 + \phi(a_j)1_{i_j}$ for $1\leq j\leq k$, we obtain
\begin{eqnarray*}
\psi(a_1a_2...a_n) &=&
\psi\big( a_1^0 ...a_{k}^0  a_{k+1}...a_n\big)\\
&+&
\sum_{1\leq j\leq k}\psi\big( a_1^0 ... 1_{i_j}a_{j+1}^0... a_k^0 a_{k+1}...a_{n}\big)\\
&+&
\sum_{1\leq j<r\leq k}
\psi\big( a_1^0 ...1_{i_{j}}a_{j+1}^0... 1_{i_r}a_{r+1}^0... a_{k}^0 a_{k+1}...a_{n}\big)\\
&+& \ldots \\
&+&
\psi\big( 1_{i_1} ... 1_{i_{k}}a_{k+1}...a_{n}\big).
\end{eqnarray*}
Now, applying (2.2) to all moments except the first one, we can reduce them to moments
of order $<n$. Therefore, by the inductive assumption, they agree with the
corresponding moments in the state $\phi$. A typical calculation is of the form
\begin{eqnarray*}
\psi\big( a_1^0 ... 1_{i_j}a_{j+1}^0 ... a_{n}\big) &=&
\psi\big( a_1^0 ... a_{j+1}^0 ... a_{n}\big)\\
&=&
\phi\big( a_1^0 ... a_{j+1}^0... a_{n}\big)\\
&=&
\phi\big( a_1^0 ... 1_{i_j}a_{j+1}^0... a_{n}\big)
\end{eqnarray*}
where the last equation is an easy property of $\phi$ for $i_{1}<i_{2}<\ldots <i_{k}$.
Thus we get
\begin{eqnarray*}
\psi(a_1a_2...a_n)&=&
\psi(a_1^0 a_2^0 ...a_{k}^0a_{k+1}...a_n)\\
&+&
\phi\big( a_1a_2... a_n\big)-
\phi\big( a_1^0 a_2^0...a_{k}^0a_{k+1}...a_n\big).
\end{eqnarray*}
However, from (2.1) we obtain
$$
\phi\big( a_1^0 a_2^0...a_{k}^0a_{k+1}...a_n\big)=0
$$
since $a_{k}^{0}\in{\rm Ker}\phi_{i_k}$ and $i_{k-1}<i_{k}>i_{k+1}$. In turn,
using (2.2) again we get
\begin{eqnarray*}
\psi\big( a_1^0 a_2^0...a_{k}^0a_{k+1}...a_n\big)&=&
\psi\big( a_1^0 a_2^0...a_{k}^0(a_{k+1}^0+1_{i_{k+1}})...a_n\big)\\
&=&
\psi\big( a_1^0 a_2^0...a_{k}^0a_{k+1}^0a_{k+2}...a_n\big)\\
&=& ...\\
&=&
\psi\big( a_1^0 a_2^0...a_n^0\big)\\
&=&
0
\end{eqnarray*}
which completes the proof of the inductive step.
This allows us to pull out the moment corresponding to the first local maximum
in the sequence $(i_1 , ... , i_n)$. In the same manner we can pull
out the moment corresponding to every local maximum.\hfill$\blacksquare$\\[10pt]

\newpage
\myownsection
\begin{center}
{\sc 3. The monotone hierarchy}
\end{center}
\indent{\par}
In this section we introduce a sequence of product states indexed by $m\in {\mathbb N}$, which gives
a discrete interpolation between the monotone product of states ($m=1$) and the free product of states
($m=\infty$).\\
\indent{\par}
{\sc Definition 3.1.}
Let $\{\mathcal{A}_{i}\}_{i\in I}$ be a family of unital *-algebras with units $\{1_{i}\}_{i\in I}$,
where $I$ is a linearly ordered set, and let $\phi: \sqcup _{i\in I}{\cal A}_{i}\rightarrow {\mathbb C}$
be a linear functional. The family $\{1_i\}_{i\in I}$ is called an {\it m-monotone family of units}
with respect to $\phi$ if $\phi(1_{i})=1$ for all $i\in I$ and
$$
\phi(a_1 \ldots a_{j-1} 1_{i_j} a_{j+1} \ldots a_n)=
\phi(a_1 \ldots a_{j-1}a_{j+1} \ldots a_n)
$$
whenever $j\leq m$ or ($m<j\leq r$ and $i_{1}\neq \ldots \neq i_{m}<\ldots <i_{r}>i_{r+1}$)
and otherwise the moment vanishes, where $a_1\in\mathcal{A}_{i_1}\cap {\rm Ker}\ \phi ,...,
a_{j-1}\in\mathcal{A}_{i_{j-1}}\cap {\rm Ker}\phi$ and $a_{j+1}\in \mathcal{A}_{i_{j+1}}, \ldots ,
a_{n} \in \mathcal{A}_{i_{n}}$.\\
\indent{\par}
Let us observe that if $\{1_{i}\}_{i\in I}$ is an $m$-monotone family of units w.r.t.
$\phi$, then each $1_{i}$ behaves as a unit when it appears at
the first $m$ positions in a given moment or if its position precedes the first
disorder in the sequence $(i_{m},\ldots , i_{n})$. Otherwise, $1_{i}$ acts as a null projection.\\
\indent{\par}
{\sc Definition 3.2.}
Let $m\in N$ and let $\{({\cal A}_{i}, \phi_{i})\}_{i\in I}$ be a family of noncommutative probability spaces,
where $I$ is a linearly ordered index set. By the $m$-{\it monotone product of states} $\phi_{i}$, $i\in I$,
we understand the linear functional $\phi:\sqcup_{i\in I}{\cal A}_{i}\rightarrow {\mathbb C}$ such that
\begin{center}
$\phi(a_1a_2...a_n)=0$,\ \ \ whenever\ \ \ $a_k\in{\cal A}_{i_k}\cap{\rm Ker}\,\phi_{i_k}$
\end{center}
and $\{1_{i}\}_{i\in I}$ is an $m$-monotone family of units w.r.t. $\phi$.
The sequence of $m$-monotone product states will be called the {\it monotone hierarchy
(of product states)}.\\
\indent{\par}
With the $m$-monotone product of states we can associate a notion of `independence' in a natural way
for every $m$.\\
\indent{\par}
{\sc Definition 3.3.}
Let $\{\mathcal{A}_{i}\}_{i\in I}$ be a family of *-subalgebras of a noncommutative probability space
$(\mathcal{A},\phi)$. A family of projections $\{1_{i}\in \mathcal{A}_{i}\}_{i\in I}$ is an
$m$-{\it monotone family of units w.r.t.} $\phi$ if $1_{i}$ is an internal unit in
${\cal A}_{i}$ for every $i\in I$ and conditions of Definition 3.1 are satisfied.
If such a family exists, then we say that $\{\mathcal{A}_{i}\}$ are {\it $m$-monotone} with respect to
$\phi$ if and only if
\begin{center}
$\phi(a_1a_2...a_n)=0$\ \ \ whenever\ \ \ $a_k\in{\cal A}_{i_k}\cap{\rm Ker}\,\phi$
\end{center}
where $k=1, \ldots , n$ and $i_{1}\ne ...\ne i_{n}$.
The associated sequence of independences will be called the {\it monotone hierarchy (of freeness)}.\\
\indent{\par}
In particular, *-subalgebras are $1$-monotone w.r.t. $\phi$ if and only if they are
monotone independent w.r.t. $\phi$. Also, if we extend Definitions 3.1 and 3.3 to include the case
$m=\infty$, we can say that *-subalgebras are $\infty$-monotone w.r.t. $\phi$ if and only if
they are free w.r.t. $\phi$. In turn, the $m$-monotone hierarchy contains models which
give examples of `freeness without identification of units'. \\
\indent{\par}
{\it Example 3.1.}
Let ${\cal H}$ be a Hilbert space, $\Omega \in {\cal H}$ - a unit vector and
$\varphi(.)=\langle . \Omega, \Omega \rangle$ - the associated state.
Then $1_{1}=1\otimes P_{\Omega}$, $1_{2}=1\otimes 1$
give a $1$-monotone family of units w.r.t. the state $\varphi \otimes \varphi$, with respect
to which ${\cal A}_{1}=B({\cal H})\otimes P_{\Omega}$ and
${\cal A}_{2}=1 \otimes B({\cal H})$ are $1$-monotone.\\
\indent{\par}
{\sc Lemma 3.4.}
{\it Let $\{\mathcal{A}_{i}\}_{i\in I}$ be a family of *-subalgebras of a noncommutative probability space
$(\mathcal{A},\phi)$ which are {\it $m$-monotone} with respect to $\phi$.
Then}
$$
\phi(a_1^0\ldots a_r^0a_{r+1}...a_n)=0
$$
{\it whenever $i_1\neq \ldots \neq i_m<i_{m+1}<\ldots <i_{r}>i_{r+1}$, where
$a_{k}^0\in {\cal A}_{i_k}\cap {\rm Ker}\phi$ for $1\leq k \leq m$ and
$a_{r}\in {\cal A}_{i_r}$ for $m+1\leq r \leq n$.}\\
\indent{\par}
{\it Proof.}
Using the definition of the $m$-monotone family of units, we obtain
\begin{eqnarray*}
\phi(a_1^0\ldots a_r^0a_{r+1}\ldots a_{n})
&=&
\phi(a_1^0\ldots a_r^0a_{r+1}^0 \ldots a_{n}) +
\phi(a_{r+1})\phi(a_1^0\ldots a_r^01_{k_{r+1}}a_{r+2}\ldots a_{n})\\
&=&
\phi(a_1^0\ldots a_r^0a_{r+1}^0a_{r+2}\ldots a_{n})
\end{eqnarray*}
We repeat the same for $a_{r+2},\ldots , a_{n}$, which leads to
$$
\phi(a_1^0\ldots a_r^0a_{r+1}\ldots a_{n})=
\phi(a_1^0\ldots a_r^0...a_n^0)
$$
and this vanishes by the definition of $m$-monotone *-subalgebras. \hfill$\blacksquare$\\
\indent{\par}
{\sc Lemma 3.5.}
{\it Let $\{\mathcal{A}_{i}\}_{i\in I}$ be a family of *-subalgebras of a noncommutative probability space
$(\mathcal{A},\phi)$ which are {\it $m$-monotone} with respect to $\phi$.
Then}
$$
\phi(a_1\ldots a_{j-1}1_{i_j}a_{j+1}...a_n)=
\phi(a_1\ldots a_{j-1}a_{j+1}\ldots a_n)
$$
{\it whenever $j=1, \ldots ,r$ and
$i_1\neq \ldots \neq i_{m}<i_{m+1}<\ldots <i_{r}>i_{r+1}$, where $a_{k}\in {\cal A}_{i_k}$ for
$1\leq k \leq n$.}\\
\indent{\par}
{\it Proof.}
For $n=1$ the statement is obvious. Suppose it holds for $1\leq s \leq n-1$. Then,
decomposing $a_{k}=a_{k}^0 +\phi(a_{k})1_{i_k}$ succesively for $k=1, \ldots , j-1$, we obtain
\begin{eqnarray*}
\phi(a_1\ldots a_{j-1}1_{i_j}a_{j+1}\ldots a_{n})
&=& \sum_{k=1}^{j-1}\phi(a_{k})\phi(a_{1}^0 \ldots a_{k-1}^0 a_{k+1}\ldots a_{j-1}1_{i_j}a_{j+1}\ldots a_{n})\\
&& +\phi(a_1^0\ldots a_{j-1}^{0}1_{i_j}a_{j+1}\ldots a_n)
\end{eqnarray*}
which, by the inductive assumption and the definition of the $m$-monotone family of units, gives
$$
\sum_{k=1}^{j-1}\phi(a_{k})\phi(a_{1}^0 \ldots a_{k-1}^0 a_{k+1}\ldots a_{j-1}a_{j+1}\ldots a_{n})
+\phi(a_1^0\ldots a_{j-1}^{0}a_{j+1}\ldots a_n)
$$
$$
=\phi(a_{1}\ldots a_{j-1}a_{j+1}\ldots a_{n})
$$
and that completes the induction proof.\hfill$\blacksquare$.\\
\indent{\par}
{\it Example 3.2.}
Suppose that $\mathcal{A}_1,\mathcal{A}_2$ are 2-monotone *-subalgebras of
$(\mathcal{A},\phi)$ and $\{1_1,1_2\}$ is a 2-monotone family of units.
Let $a\in\mathcal{A}_1,\,b\in\mathcal{A}_2$. It can be easily shown that mixed
moments of orders $\leq 4$ agree with the corresponding moments in the free case.
Namely
\begin{eqnarray*}
\phi(ab)&=&\phi(a)\phi(b),\\
\phi(aba)&=& \phi(a^{2})\phi(b),\\
\phi(abab)&=&\phi(a^{2})\phi^{2}(b)+\phi(b^{2})\phi^{2}(a)-\phi^{2}(a)\phi^{2}(b).
\end{eqnarray*}
Let us compute the `alternating' mixed moment of order 5.
We have
\begin{eqnarray*}
\phi(ababa)&=& \phi(a^0baba)+\phi(a)\phi(baba)\\
&=& \phi(a^0b^0aba)+\phi(b)\phi(a^0aba)+\phi(a)\phi(b^0aba)+\phi(a)\phi(b)\phi(aba)\,.
\end{eqnarray*}
Let us observe that the first term in this expression vanishes
by Lemma 3.4. Besides, easy computations of lower-order moments give
\begin{eqnarray*}
\phi(a^0aba) &=& \phi(b)\phi(a^3)-\phi(b)\phi(a)\phi(a^2),\\
\phi(b^0aba) &=& \phi^2(a)\phi(b^2)-\phi^2(a)\phi^2(b),\\
\phi(aba) &=& \phi(b)\phi(a^2).
\end{eqnarray*}
Thus we finally obtain
$$
\phi(ababa) = \phi^3(a)\phi(b^2) - \phi^3(a)\phi^2(b) + \phi(a^3)\phi^2(b),
$$
which gives the lowest order moment for $m=2$
which differs from the corresponding moments in both the free and the monotone cases. \\[10pt]

\myownsection
\begin{center}
{\sc 4. Hilbert space representations}
\end{center}
We begin this section with recalling the definition of
the free product of noncommutative probability spaces [V].
Let $I$ be a linearly ordered set of indices and let
$(\mathcal{A}_i,\phi_i)_{i\in I}$ be a family of
$\mathcal{C}^*$-noncommutative probability spaces.
For each $i\in I$ we denote by $(\mathcal{H}_i,\pi_i,\xi_i)$ the GNS triple
of $(\mathcal{A}_i, \phi_{i})$, i.e. $\mathcal{H}_i$
is a Hilbert space, $\xi_i$ is a cyclic (unit) vector in
$\mathcal{H}_i$ and $\pi_i:\mathcal{A}_i\rightarrow\mathcal{B}(\mathcal{H}_i)$ is a *-homomorphism, such that
$$
\phi_i(a)=\langle \pi_i(a)\xi_i , \xi_i \rangle\
$$
for all $a\in\mathcal{A}_i$, where $\langle\cdot\,,\cdot\rangle$ denotes the scalar product
in $\mathcal{H}_i$ (for simplicity, we use the same notation for all scalar products).

Let $\mathcal{H}_i^0$ be the orthogonal complement in $\mathcal{H}_i$ of the one-dimensional
subspace generated by $\xi_i$, i.e. $\mathcal{H}_i^0=\mathcal{H}_i \ominus \mathbb{C} \omega_i$.
For any $h\in\mathcal{H}_i$, let $h^0$ denote the orthogonal projection of $h$ on $\mathcal{H}_i^0$.
Besides, let $\mathcal{H}$ denote the Hilbert space
\begin{equation}
\mathcal{H}=\mathbb{C}\,\Omega\oplus\bigoplus_{n=1}^{\infty}\bigoplus_{i_1 \ne i_2 \ne ... \ne i_n}
\mathcal{H}_{i_1}^0\otimes\mathcal{H}_{i_2}^0\otimes ... \otimes\mathcal{H}_{i_n}^0
\end{equation}
with the canonical scalar product, where $\Omega$ is the so-called vacuum (unit) vector.
The pair $(\mathcal{H},\Omega)$ is called the {\it free product} of $(\mathcal{H}_i,\xi_i)$.

On ${\cal H}$ we define a *-representation $\lambda_i:\mathcal{A}_i\rightarrow\mathcal{B}(\mathcal{H})$
of each algebra $\mathcal{A}_i$, $i\in I$ as follows:
\begin{eqnarray*}
\lambda_i(a)\Omega&=&(\pi_i(a)\xi_i)^0+\langle \pi_i(a)\xi_i,\xi_i \rangle\Omega\\
\lambda_i(a)(h_1\otimes...\otimes h_n)&=&
\begin{cases}
(\pi_i(a)\xi_i)^0\otimes h_1\otimes...\otimes h_n +
\langle\pi_i(a)\xi_i,\xi_i\rangle h_1\otimes...\otimes h_n & \\
\ \ \ \ \ \ \ \ \ \ \ \ \ \ \ \ \ \ \ \ \ \ \ \ \ \ \ \ \ \ \ \ \ \ \ \ \ \ \ \ \ \ \ \ \ \ \ \ \
\hfill{\footnotesize{\text{if $i\ne i_1$}}}\\
(\pi_i(a)h_1)^0\otimes h_2\otimes...\otimes h_n + \langle\pi_i(a)h_1,\xi_i\rangle h_2\otimes...\otimes h_n & \\
\ \ \ \ \ \ \ \ \ \ \ \ \ \ \ \ \ \ \ \ \ \ \ \ \ \ \ \ \ \ \ \ \ \ \ \ \ \ \ \ \ \ \ \ \ \ \ \ \
\hfill{\footnotesize{\text{if $i=i_1$}}}\\
\end{cases}
\end{eqnarray*}
Finally, on $\mathcal{B}(\mathcal{H})$ we define the so-called {\it vacuum state}
$\varphi(\cdot)=\langle \cdot\Omega,\Omega\rangle$.

In order to introduce the $m$-monotone product of Hilbert spaces, we want to characterize
admissible simple tensors which appear in the product.
For that purpose we introduce the associated families of sequences $I_{n}(m)\subset I^{n}$.
We set
$$
I_n(m)=\{(i_1,...,i_n); i_1 \ne i_2 \ne ... \ne i_n\} \ \ \text{for $n\leqslant m$},
$$
whereas
$$
I_n(m)=\{(i_1,...,i_n); i_1 > ... > i_{n-m+1} \ne i_{n-m+2} \ne ... \ne i_n\}\ \ \text{for $n>m$}.
$$
Thus, in contrast to the free product case, all sequences $(i_{1},i_{2}, \ldots , i_{n})$
subject to $i_{1}\neq i_{2}\neq \ldots \neq i_{n}$ are allowed only if $n\leq m$.
In turn, if $n>m$, then the first $n-m+1$ indices must be in a decreasing order.\\
\indent{\par}
{\sc Definition 4.1.}
Let $\mathcal{H}^{(m)}$ be the closed subspace of the Hilbert space $\mathcal{H}$ of the form
\begin{equation}
\mathcal{H}^{(m)}=\mathbb{C}\,\Omega\oplus\bigoplus_{n=1}^{\infty}
\bigoplus_{(i_1,...,i_n)\in I_n(m)} \mathcal{H}_{i_1}^0\otimes\mathcal{H}_{i_2}^0\otimes ...
\otimes\mathcal{H}_{i_n}^0
\end{equation}
The pair $(\mathcal{H}^{(m)},\Omega)$ will be called the {\it $m$-monotone product of Hilbert spaces}
$(\mathcal{H}_i,\xi_i)$.\\
\indent{\par}
It is natural to compare this truncation with the $m$-{\it free product of Hilbert spaces} [F-L] given by
$$
{\cal H}_{(m)}={\mathbb C}\Omega \oplus \bigoplus _{n=1}^{m} \bigoplus_{i_{1}\neq i_{2}\neq \ldots \neq i_{n}}
{\cal H}_{i_1}^{0}\otimes {\cal H}_{i_2}^{0}\otimes \ldots \otimes {\cal H}_{i_n}^{0}
$$
and observe that, in contrast to ${\cal H}^{(m)}$, ${\cal H}_{(m)}$ is a truncation of the
free product of Hilbert spaces which takes into account only sequences
$(i_{1}, i_{2}, \ldots , i_{n})$ of lenght $\leq m$. In the case of ${\cal H}^{(m)}$,
the truncation is more delicate and depends on the order in which indices appear. Thus, for instance, if
$I$ is infinite, simple tensors of arbitrary lenght appear in ${\cal H}^{(m)}$.\\
\indent{\par}
{\it Example 4.1.}
Let $I=\{1,2\}$. If $m=1$, we obtain
$$
{\cal H}^{(1)}={\mathbb C}\Omega \oplus {\cal H}_{1}^{0} \oplus {\cal H}_{2}^{0}\oplus
({\cal H}_{2}^{0}\otimes {\cal H}_{1}^{0})
$$
thus the $1$-monotone product of ${\cal H}_{1}$ and ${\cal H}_{2}$ coincides with the monotone Fock space.
In turn, the $2$-monotone product of ${\cal H}_{1}$ and ${\cal H}_{2}$ is of the form
$$
{\cal H}^{(2)}={\mathbb C}\Omega \oplus {\cal H}_{1}^{0} \oplus {\cal H}_{2}^{0}\oplus
({\cal H}_{2}^{0}\otimes {\cal H}_{1}^{0})\oplus ({\cal H}_{1}^{0}\otimes {\cal H}_{2}^{0})
\oplus
({\cal H}_{2}^{0}\otimes {\cal H}_{1}^{0}\otimes {\cal H}_{2}^{0})
$$
thus it is a direct sum of of all simple tensors of order $\leq 2$ and only
one simple tensor of order $3$.\\
\indent{\par}
Let us introduce the closed subspace of $\mathcal{H}^{(m)}$ given by
$$
\mathcal{H}^{(m)}(i)={\cal H}_{(m-1)}
\oplus\Bigg(\bigoplus_{n=m}^{\infty}
\bigoplus_{\stackrel{(i_1,...,i_n)\in I_n(m)}{\tiny i_1\,\leqslant\ i}}
\mathcal{H}_{i_1}^0\otimes...\otimes\mathcal{H}_{i_n}^0\Bigg),
$$
where $i\in I$, which allows us to define suitable truncations of the free product representation.
Let representations $\lambda_i^{(m)}$ of $\mathcal{A}_i$ on $\mathcal{H}^{(m)}$ be defined by
$$
\lambda_i^{(m)}:\mathcal{A}_i \rightarrow \mathcal{B}(\mathcal{H}^{(m)})\ \ \ ,\ \ \ \ \ \ \ \ \lambda^{(m)}_i(a)h=
\begin{cases}
\lambda_i(a)h & \text{if $h\in \mathcal{H}^{(m)}(i)$}\\
\ \ \ 0 & \text{if  $h\in \mathcal{H}^{(m)}(i)^{\perp}$}
\end{cases}.
$$
Thus, $\lambda_i^{(m)}$ is a truncation of $\lambda_i$ for every $m$. The `level' of this truncation
depends on the order of Hilbert spaces at all sites except the last $m$ sites where the order is irrelevant
as in the case of representations $\lambda_i$. Of course, the free product $*_{i\in I}\lambda_i^{(m)}$
is then a truncation of the free product representation $*_{i\in I}\lambda_i$ and agrees with the latter
on the $m$-free product of Hilbert spaces.\\
\indent{\par}
{\sc Theorem 4.2.}
{\it Algebras $\{\lambda_i^{(m)}(\mathcal{A}_i)\}_{i\in I}$ are m-monotone *-subalgebras
of the noncommutative probability space $(\mathcal{B}(\mathcal{H}^{(m)}), \varphi)$,
where $\varphi$ is the vacuum state on $\mathcal{B}(\mathcal{H}^{(m)})$.}\\
\indent{\par}
{\it Proof.}
Let us fix $m\in \mathbb{N}$. We will show that $\{\lambda_{i}^{(m)}(1_i)\}_{i\in I}$ is an
$m$-monotone family of units for $\{\lambda_{i}^{(m)}(\mathcal{A}_i)\}_{i\in I}$. Notice that
$$
\lambda_{i}^{(m)}(1_i)h=
\begin{cases}
h & \text{if $h\in \mathcal{H}^{(m)}(i)$}\\
0 & \text{if  $h\in \mathcal{H}^{(m)}(i)^{\perp}$}
\end{cases}.
$$
Besides, let $a_1\in\mathcal{A}_{i_1},...,a_n\in\mathcal{A}_{i_n}$ be such that
$i_1\ne i_2\ne...\ne i_n$. For notational simplicity we denote $X_{k}=\lambda_{i_k}^{(m)}(a_{i_k})$.
Then
\begin{eqnarray*}
\varphi(X_1 ... X_{j-1}\,\lambda_{i}^{(m)}(1_{i_j})\,X_{j+1} ... X_n)
&=& \langle X_1 ... X_{j-1}\,\lambda_{i}^{(m)}(1_{i_j})\,X_{j+1} ... X_n\Omega,\Omega \rangle\\
&=& \langle \Omega,X_n^* ... X_{j+1}^*\,\lambda_{i}^{(m)}(1_{i_j})\,X_{j-1}^* ... X_1^* \Omega \rangle\,.
\end{eqnarray*}
Let $X_{k}\in {\rm Ker}\phi$ for $k=1, \ldots , j-1$. If $(i_j,...,i_1)\in I_j(m)$,
then it can be seen that $X_{j-1}^* ... X_1^*\Omega$ is a simple tensor from $\mathcal{H}^{(m)}(i_j)$,
hence
\begin{eqnarray*}
\varphi(X_1 ... X_{j-1}\,\lambda_{i}^{(m)}(1_{i_j})\,X_{j+1} ... X_n) &=&
\langle \Omega,X_n^* ... X_{j+1}^*\,\lambda_{i}^{(m)}(1_{i_j})\,X_{j-1}^* ... X_1^*\Omega \rangle\\
&=& \langle\Omega,X_n^* ... X_{j+1}^*X_{j-1}^*...X_1^*\Omega \rangle\\
&=& \varphi(X_1 ... X_{j-1}X_{j+1} ... X_n)\,.
\end{eqnarray*}
If, in turn, $(i_j,...,i_1)\notin I_j(m)$, then
$X_{j-1}^* ... X_1^*\Omega \in \mathcal{H}^{(m)}(i_j)^{\perp}$, hence
$$
\lambda_{i}^{(m)}(1_{i_j})\,X_{j-1}^* ... X_1^*\Omega=0,
$$
which proves that $\{\lambda_{i}^{(m)}(1_i)\}_{i\in I}$ is an $m$-monotone family of units.

Let us show now that variables $X_1,...,X_n$ satisfy conditions of Definition 3.3.
Assume that $X_1,...,X_n\in {\rm Ker}\,\phi$. Then
$X_n^*...X_1^*\Omega\perp \mathbb{C}\Omega$ since it is a simple tensor from ${\cal H}^{(m)}(i_n)$,
and thus
$$
\varphi(X_1...X_n)=\langle\,\Omega,X_n^* ... X_1^*\Omega\,\rangle=0\,,
$$
which completes the proof.\hfill{$\blacksquare$}\\
\indent{\par}
{\sc Corollary 4.3.}
{\it The $m$-monotone product of states is a state which agrees with the free product of states
on words of lenght less than or equal to $2m$.}\\
\indent{\par}
{\it Proof.}
It follows from the GNS representation in the proof of Theorem 4.2
that the $m$-monotone product of states is positive, thus it is a state.
Moroever, by definition, it satisfies the freeness condition and the units act as an identified unit
at the first $m$ places in the moment.
From this and the fact that any state is a hermitian functional it follows that
units act as an identified unit at the last $m$ places in the moment. Therefore,
the moments have to agree if their lenghts are smaller or equal to $2m$. \hfill$\blacksquare$\\[10pt]
\myownsection
\begin{center}
{\sc 5. Combinatorics}
\end{center}
This section is devoted to the combinatorics of the monotone hierarchy.
It is based on ordered non-crossing partitions with disorders starting
at different depths (compare with the combinatorics of non-crossing
parititions in free probability [S]).

An {\it ordered partition} $P$ of the set $\{1,2, \dots , n\}$ is a tuple
$(P_{1},P_{2}, \ldots , P_{k})$ of disjoint
non-empty subsets  called {\it blocks} such that
$P_{1}\cup P_{2}\cup \ldots \cup P_{k}=\{1,2,\ldots , n\}$. For any block $P_{r}$ the number
$r$ is the {\it order} of $P_{r}$, which can also be interpreted as {\it color}.
If blocks are two-element sets, then $P$ is called an {\it ordered pair-partition}.
An ordered partition $S$ is called {\it crossing} if and only if
there exist $s,s'\in P_{i}$ and $r,r'\in P_{j}$ for some $i\neq j$ such that
$s<r<s'<r'$. If $P$ is not a crossing partition, then
it is called {\it non-crossing}. The collection of ordered non-crossing partitions
(pair-partitions) of the set $\{1,2, \ldots , n\}$ will be denoted ${\cal ONC}_{n}$ (${\cal ONC}^2_n$).
Block $P_i$ is called {\it inner} with respect
to block $P_{j}$, which we denote $P_j<P_i$, if and only if there exist two
numbers $p,q\in P_j$ such that the interval $(p,q)\cap P_{j}=\emptyset$ and
for every $s\in P_{i}$ we have $p<s<q$. Then $P_{j}$ will
be called {\it outer} with respect to $P_{i}$. By the {\it depth} of block $P_i$,
denoted $d(P_i)$, we will understand the number of blocks of $P$ which are outer w.r.t.
$P_i$ or coincide with $P_i$. In particular, if there are no blocks which are outer w.r.t. $P_{i}$,
then $d(P_i)=1$. \\
\indent{\par}
{\sc Definition 5.1.}
By ${\cal ONC}_{n}(m)$ we will denote the collection of all {\it ordered non-crossing
partitions} $P=(P_1,P_2,...,P_k)$ of the set $\{1,2,...,n\}$ such that for all
$i,j\in \{1,...,k\}$ the following implication holds:
$$
d(P_j)\geq m\;\;{\rm and}\;\; P_j < P_i\; \Longrightarrow\; j < i .
$$
In other words, the color of blocks is a monotone function of their depths starting from
depth $m$.
In particular, we call ${\cal ONC}_{n}(1)$ {\it monotone non-crossing partitions}
(we also denote them ${\cal MNC}_{n}$).
The corresponding pair partitions will be denoted ${\cal ONC}_{n}^{2}(m)$
and ${\cal MNC}_{n}^{2}$, respectively.\\
\indent{\par}
{\it Example 5.1.} For instance, for $n=8$ the tuples
\begin{center}
$P$=(\{4,7\},\{1,8\},\{2,3\},\{5,6\}) and $R$=(\{1,8\},\{2,3\},\{4,7\},\{5,6\})
\end{center}
are different ordered pair-partitions of the set \{1,2,...,8\}.
The corresponding diagrams are drawn in Figure 1, where
blocks are labelled by colors associated with their orders.\\

\unitlength=1mm
\special{em:linewidth 0.4pt}
\linethickness{0.4pt}
\begin{picture}(120.00,50.00)(-15.00,5.00)

\put(0.00,31.00){\line(0,1){18.00}}
\put(0.00,49.00){\line(1,0){42.00}}
\put(42.00,31.00){\line(0,1){18.00}}
\put(20.00,49.50){\small{\it{2}}}

\put(6.00,31.00){\line(0,1){6.00}}
\put(6.00,37.00){\line(1,0){6.00}}
\put(12.00,31.00){\line(0,1){6.00}}
\put(8.50,37.50){\small{\it{3}}}

\put(18.00,31.00){\line(0,1){12.00}}
\put(18.00,43.00){\line(1,0){18.00}}
\put(36.00,31.00){\line(0,1){12.00}}
\put(26.00,44.00){\small{\it{1}}}

\put(24.00,31.00){\line(0,1){6.00}}
\put(24.00,37.00){\line(1,0){6.00}}
\put(30.00,31.00){\line(0,1){6.00}}
\put(26.00,37.50){\small{\it{4}}}

\put(0.00,31.00){\circle*{1.00}}
\put(6.00,31.00){\circle*{1.00}}
\put(12.00,31.00){\circle*{1.00}}
\put(18.00,31.00){\circle*{1.00}}
\put(24.00,31.00){\circle*{1.00}}
\put(30.00,31.00){\circle*{1.00}}
\put(36.00,31.00){\circle*{1.00}}
\put(42.00,31.00){\circle*{1.00}}

\put(-0.90,27.00) {\footnotesize 1}
\put(05.10,27.00) {\footnotesize 2}
\put(11.10,27.00) {\footnotesize 3}
\put(17.10,27.00) {\footnotesize 4}
\put(23.10,27.00) {\footnotesize 5}
\put(29.10,27.00) {\footnotesize 6}
\put(35.10,27.00) {\footnotesize 7}
\put(41.10,27.00) {\footnotesize 8}

\put(60.00,31.00){\line(0,1){18.00}}
\put(60.00,49.00){\line(1,0){42.00}}
\put(102.00,31.00){\line(0,1){18.00}}
\put(80.00,49.50){\small{\it{1}}}

\put(66.00,31.00){\line(0,1){6.00}}
\put(66.00,37.00){\line(1,0){6.00}}
\put(72.00,31.00){\line(0,1){6.00}}
\put(68.50,37.50){\small{\it{2}}}

\put(78.00,31.00){\line(0,1){12.00}}
\put(78.00,43.00){\line(1,0){18.00}}
\put(96.00,31.00){\line(0,1){12.00}}
\put(86.00,44.00){\small{\it{3}}}

\put(84.00,31.00){\line(0,1){6.00}}
\put(84.00,37.00){\line(1,0){6.00}}
\put(90.00,31.00){\line(0,1){6.00}}
\put(86.00,37.50){\small{\it{4}}}

\put(60.00,31.00){\circle*{1.00}}
\put(66.00,31.00){\circle*{1.00}}
\put(72.00,31.00){\circle*{1.00}}
\put(78.00,31.00){\circle*{1.00}}
\put(84.00,31.00){\circle*{1.00}}
\put(90.00,31.00){\circle*{1.00}}
\put(96.00,31.00){\circle*{1.00}}
\put(102.00,31.00){\circle*{1.00}}

\put(59.10,27.00) {\footnotesize 1}
\put(65.10,27.00) {\footnotesize 2}
\put(71.10,27.00) {\footnotesize 3}
\put(77.10,27.00) {\footnotesize 4}
\put(83.10,27.00) {\footnotesize 5}
\put(89.10,27.00) {\footnotesize 6}
\put(95.10,27.00) {\footnotesize 7}
\put(101.10,27.00) {\footnotesize 8}

\put(60,19){\small{$R$=(\{1,8\},\{2,3\},\{4,7\},\{5,6\})}}
\put(-2,19){\small{$P$=(\{4,7\},\{1,8\},\{2,3\},\{5,6\})}}
\put(16,5){{\it Figure 1.} Diagrams of partitions $P$ and $R$.}
\end{picture}
\\
\indent{\par}
In the case of partitions $P$ and $R$ in Figure 1 all blocks but $\{1,8\}$ are inner w.r.t. $\{1,8\}$, but
$\{5,6\}$ is the only block which is inner w.r.t. $\{4,7\}$.
The depths are given by
$d(\{1,8\})=1$, $d(\{2,3\})=d(\{4,7\})=2$, $d(\{5,6\})=3$.
Note that $R\in {\cal MNC}_{8}^{2}$ since its coloring is monotone, whereas
$P\in {\cal ONC}_{8}^{2}(2)\setminus {\cal MNC}_{8}^{2}$ since the coloring of the whole partition
is not monotone, however it is monotone starting from depth $2$.\\
\indent{\par}
{\sc Proposition 5.2.}
{\it For every $k\in {\mathbb N}$ we have
$\left|{\cal MNC}^2_{2k}\right|\ =\ (2k-1)!!$.}\\
\indent{\par}
{\it Proof.}
The proof is elementary and is left to the reader.\hfill{$\blacksquare$}\\
\indent{\par}
{\sc Definition 5.3.}
Let $i_1,...,i_n\in\mathbb{N}$ and set $\{i_{1},\ldots , i_{n}\}=
\{k_{1}, \ldots , k_{r}\}$, where $k_{1}<k_{2}<\ldots <k_{r}$.
By an ordered partition {\it associated} with the tuple $(i_{1},\ldots , i_{n})$
we understand the partition $P=(P_{1},\ldots , P_{r})$ given by
$P_{j}=\{s:i_{s}=k_{j}\}$. We then write $P\sim (i_{1}, \ldots , i_{n})$.\\
\indent{\par}
{\it Example 5.2.}
Let $i_{1}=i_{4}=2$, $i_{2}=i_{5}=4$, $i_{3}=1$. Then $(P_{1},P_{2},P_{3}) \sim (2,4,1,2,4)$,
where $P_{1}=\{3\}$, $P_{2}=\{1,4\}$, $P_{3}=\{2,5\}$.\\
\indent{\par}
We shall also need a 'continuous' version of Definition 5.3 with
the tuple of indices replaced by a tuple of characteristic functions. For that purpose, on the set
\begin{equation}
\Theta =\{\chi_{(s,t]}: 0\leq s <t <\infty\}
\end{equation}
introduce a partial order by setting $f=\chi_{(s_1,t_1]} < \chi_{(s_2,t_2]}=g$ if and only if $t_1\leq s_2$
and $f\leq g$ iff $f<g$ or $f=g$.\\
\indent{\par}
{\sc Definition 5.4.}
Let $f_{1}, f_{2}, \ldots , f_{n}\in \Theta$ have identical or disjoint supports.
We say that $P\in {\cal OP}_{n}$ is {\it compatible} with the tuple
$(f_{1},f_{1}, \ldots , f_{n})$ if and only if
\indent{\par}1. $i,j\in P_{k}\Longrightarrow f_{i}=f_{j}$,
\indent{\par}2. $i\in P_{k}, j\in P_{l}$ and $k<l\;\Longrightarrow \; f_{i}\leq f_{j}$.\\
We then write $P\sim (f_{1},f_{2}, \ldots , f_{n})$.\\
\indent{\par}
{\it Example 5.3.}
Let $f_{1}=f_{2}=f_{5}=f_{6}=\chi_{(0,1]}$ and $f_{3}=f_{4}=\chi_{(1,2]}$. Then the partition
$P=(\{1,6\},\{2,5\}, \{3,4\})$ is compatible with the tuple $(f_{1},f_{2}, \ldots , f_{6})$,
whereas $P'=(\{3,4\}, \{1,6\}, \{2,5\})$ is not since, for instance, block $\{3,4\}$ precedes
block $\{1,6\}$ and $f_{3}=f_{4}>f_{1}=f_{6}$.\\
\indent{\par}
The combinatorial formulas derived in this paper can be expressed in terms of
the usual non-crossing partitions.
If $P=(P_{1},P_{2}, \ldots , P_{k})\in {\cal ONC}_{n}$ and
$\pi=\{\pi_{1}, \ldots , \pi_{k}\}\in {\cal NC}_{n}$ have the same blocks, then we shall
write $\pi\sim P$ (by abuse of notation).

Let us assume now that we have a pair of objects: a partition
$\pi= \{\pi_{1},\pi_{2} \ldots , \pi_{k}\}\in {\cal NC}_{2k}^{2}$ and
a tuple of functions $(f_{1}, f_{2}, \ldots , f_{2k})$, where $f_{i}=\xi_{(s_{i},t_{i}]}$,
$i=1, \ldots , 2k$, have identical or disjoint supports and such that
$\pi \sim P\sim (f_{1}, f_{2}, \ldots , f_{2k})$ for some ordered partition $P$.
This implies that if $\{p,q\}=\pi_{i}$
for some $i$, then $f_{p}=f_{q}:=f^{(i)}=\chi_{(s^{(i)},t^{(i)}]}$.
Thus $f^{(1)},f^{(2)},  \ldots , f^{(k)}$ will denote characteristic functions associated with blocks
$\pi_{1},\pi_{2} \ldots , \pi_{k}$ and $s^{(1)},s^{(2)}, \ldots , s^{(k)}$ as well as
$t^{(1)},t^{(2)}, \ldots , t^{(k)}$ will be the left and the right endpoints of their supports,
respectively.

Computations will be based on counting inner blocks which are 'compatible'
with the supports of associated characteristic functions.
For that purpose, to every block $\pi_{i}$ of $\pi$ we assign the number of its blocks which are inner w.r.t.
$\pi_{i}$ and are associated with functions of the same supports as $f^{(i)}$, namely
\begin{equation}
{\rm Inn}(\pi_{i}) = \#(\pi_{j}: \pi_{j}>\pi_{i}\; {\rm and} \;f^{(j)}=f^{(i)}).
\end{equation}
Note that ${\rm Inn}(\pi_{i})$ depends on the tuple $(f_{1}, f_{2}, \ldots , f_{2k})$, but this fact
is supressed in the notation.

Finally, to every $\pi$ and $(f_{1},f_{2} \ldots , f_{2k})$ we associate the number
\begin{equation}
c_{\pi}(f_{1}, f_{2}, \ldots , f_{2k})=\#(P\in {\cal MNC}_{2k}^{2}: \pi\sim P\sim (f_{1},f_{2}, \ldots , f_{2k}))
\end{equation}
which gives the number of `admissible' colorings of $\pi$.\\
\indent{\par}
{\it Example 5.4.}
Let $\pi\in {\cal NC}_{8}^{2}$ be given by the diagram in Figure 2 with supports $f=\chi_{(0,t]}$,
and $g=\chi_{(t,t']}$ associated with blocks as shown. Set $\pi_{1}=\{1,8\}$,
$\pi_{2}=\{2,3\}$, $\pi_{3}=\{4,7\}$, $\pi_{4}=\{5,6\}$.\\
\unitlength=1mm
\special{em:linewidth 0.4pt}
\linethickness{0.4pt}
\begin{picture}(120.00,42.00)(-50.00,5.00)

\put(0.00,21.00){\line(0,1){18.00}}
\put(0.00,39.00){\line(1,0){42.00}}
\put(42.00,21.00){\line(0,1){18.00}}
\put(20.00,40.50){\small{\it{f}}}

\put(6.00,21.00){\line(0,1){6.00}}
\put(6.00,27.00){\line(1,0){6.00}}
\put(12.00,21.00){\line(0,1){6.00}}
\put(8.50,28.50){\small{\it{g}}}

\put(18.00,21.00){\line(0,1){12.00}}
\put(18.00,33.00){\line(1,0){18.00}}
\put(36.00,21.00){\line(0,1){12.00}}
\put(26.00,34.50){\small{\it{f}}}

\put(24.00,21.00){\line(0,1){6.00}}
\put(24.00,27.00){\line(1,0){6.00}}
\put(30.00,21.00){\line(0,1){6.00}}
\put(26.00,28.50){\small{\it{g}}}

\put(0.00,21.00){\circle*{1.00}}
\put(6.00,21.00){\circle*{1.00}}
\put(12.00,21.00){\circle*{1.00}}
\put(18.00,21.00){\circle*{1.00}}
\put(24.00,21.00){\circle*{1.00}}
\put(30.00,21.00){\circle*{1.00}}
\put(36.00,21.00){\circle*{1.00}}
\put(42.00,21.00){\circle*{1.00}}

\put(-0.90,17.00) {\footnotesize 1}
\put(05.10,17.00) {\footnotesize 2}
\put(11.10,17.00) {\footnotesize 3}
\put(17.10,17.00) {\footnotesize 4}
\put(23.10,17.00) {\footnotesize 5}
\put(29.10,17.00) {\footnotesize 6}
\put(35.10,17.00) {\footnotesize 7}
\put(41.10,17.00) {\footnotesize 8}

\put(-15,5){{\it Figure 2.} Diagram of the partition $\pi$.}
\end{picture}
\\[10pt]
Block $\pi_{1}$ has one inner block with support equal to $f$, thus ${\rm Inn}(\pi_{1})=1$,
whereas the remaining blocks have no inner blocks with identical supports, thus
${\rm Inn}(\pi_{2})={\rm Inn}(\pi_{3})={\rm Inn}(\pi_{4})=0$. Note also that
$c_{\pi}(f,g,g,f,g,g,f,f)=3$ (block $\pi_{1}$ must be colored by $1$, block $\pi_{2}$ can
be colored by $2,3$ or $4$ and for every such choice there is only one coloring
of the remaining blocks since $\pi_{3}<\pi_{4}$).\\[10pt]
\myownsection
\begin{center}
{\sc 6. Lemmas on mixed moments}
\end{center}
In the lemmas given below we assume that
$\{\mathcal{A}_{i}\}_{i\in I}$ is a family of *-subalgebras of a noncommutative probability space
$(\mathcal{A},\phi)$ which are {\it $m$-monotone} with respect to $\phi$ and that
$a_{k}\in {\cal A}_{i_k}$ for $1\leq k \leq n$. \\
\indent{\par}
{\sc Lemma 6.1.}
{\it Let $\{\mathcal{A}_i\}_{i\in I}$ be $m$-monotone *-subalgebras of
$(\mathcal{A},\phi)$ and let $a_1\in\mathcal{A}_{i_1} , ... , a_n\in\mathcal{A}_{i_n}$ be such that
$i_1\ne i_2\ne ... \ne i_n$. If there exists $j\in\{1,...,n\}$, for which
$\phi(a_j)=0$ and $i_j\ne i_k$ when $j\ne k$, then $\phi(a_1a_2...a_n)=0$.}\\
\indent{\par}
{\it Proof.}
We use induction w.r.t. $n$. For $n=1$ the statement is obvious.
Suppose that it holds for $l<n$. We have
\begin{eqnarray*}
\phi(a_1a_2...a_n)&=& \phi(a_1^0a_2...a_n)+\phi(a_1)\phi(1_{i_1}a_2...a_n)\\
&=&\phi(a_1^0a_2...a_n)+\phi(a_1)\phi(a_2...a_n)\,.
\end{eqnarray*}
If $j\ne 1$, then $\phi(a_2...a_n)=0$ from the inductive assumption.
If $j=1$, then $\phi(a_1)=0$, which gives
$$
\phi(a_1a_2...a_n)=\phi(a_1^0a_2...a_n)=\phi(a_1^0a_2^0a_3...a_n)+\phi(a_2)\phi(a_1^01_{i_2}a_3...a_n)
$$
As before, consider two cases. In the first case, namely $j=2$ we have $\phi(a_2)=0$, whereas in the second case,
for $j\ne 2$ we have two possibilities:\\

\noindent
1. $m\leqslant 1$. If $i_1 < i_2$, then $\phi(a_1^01_{i_2}a_3...a_n)= \phi(a_1^0a_3...a_n)=0$
from the inductive assumption.  In turn, if $i_1 > i_2$, then $\phi(a_1^01_{i_2}a_3...a_n)=0$, which
follows from Definition 3.1.\\

\noindent
2. $m > 1$. Then $\phi(a_1^01_{i_2}a_3...a_n)= \phi(a_1^0a_3...a_n)=0$ from the inductive assumption.\\

\noindent
Thus, in all considered cases, either $\phi(a_1^01_{i_2}a_3...a_n)=0$, or $\phi(a_2)=0$, hence
$$
\phi(a_1a_2...a_n)=\phi(a_1^0a_2^0a_3...a_n).
$$
By continuing this reasoning we obtain $\phi(a_1a_2...a_n)=\phi(a_1^0a_2^0...a_n^0)$, which
is equal to zero by Definition 3.2.\hfill{$\blacksquare$}\\
\indent{\par}
{\sc Lemma 6.2.}
{\it Suppose that $(i_{1},\ldots , i_{n})\sim P$, where $i_{1}\neq i_{2}\neq \ldots \neq i_{n}$ and
where $P$ is a crossing partition.
Then the mixed moment $\phi(a_{1}\ldots a_{n})$ is expressed
in terms of products of at least $p+1$ marginal moments, where $p$ is the number of blocks of $P$.}\\
\indent{\par}
{\it Proof.}
For $n=4$, where the induction starts, the statement is true since
$\phi(abab)=\phi(a^{2})\phi^{2}(b)$ for $m=1$ (monotone case) and
$$
\phi(abab)=\phi^{2}(b)\phi(a^{2})-\phi^{2}(b)\phi^{2}(a)+
\phi^{2}(a)\phi(b^{2})
$$
for $m>1$ (then, by Corollary 4.4, moments of order 4  agree with those in the free case).
Let $n>4$ and suppose that $i_1\neq \ldots \neq i_m<i_{m+1}<\ldots <i_r>i_{r+1}$. In the formula
$$
\phi(a_{1}a_{2}\ldots a_{n})=\phi(a_{1}^0a_{2}\ldots a_{n})+\phi(a_{1})\phi(a_{2}\ldots a_{n})
$$
the second summand is a product of at least $p+1$ factors since
either $(i_{2} ,\ldots , i_{n})$ corresponds to a non-crossing partition $P'$, and then
$P'$ must have $p$ blocks, or $(i_{2}, \ldots , i_{n})$ corresponds to a crossing partition, in which case
we use the inductive assumption. In turn, we can write
$$
\phi(a_{1}^0a_{2}\ldots a_{n})=\phi(a_{1}^0a_{2}^0a_{3}\ldots a_{n})+\phi(a_{2})\phi(a_{1}^0a_{3}\ldots a_{n})
$$
and use a similar argument to conclude that the second summand factorizes into at least $p+1$ factors.
Continuing this reasoning, we can reduce the proof to the case of $\phi(a_{1}^0\ldots a_r^0a_{r+1}\ldots a_n)$,
which vanishes by Lemma 3.4.\hfill {$\blacksquare$}\\
\indent{\par}
{\sc Lemma 6.3.}
{\it Suppose that $(i_{1},\ldots , i_{n})\sim P\in {\cal ONC}_{n}(m)$, where $i_{1}\neq i_{2}\neq \ldots \neq i_{n}$.
Then the mixed moment $\phi(a_{1}\ldots a_{n})$ is a product of exactly
$p$ marginal moments, where $p$ is the number of blocks of $P$.}\\
\indent{\par}
{\it Proof.}
We use induction w.r.t. $n$. For $n\leq 2m$, we get the same expression for the moments as in the free case
by Corollary 4.4 and thus a product of $p$ marginal moments [Sp]. Thus let $n>2m$ and suppose that
$i_{1}\neq \ldots \neq i_{m}<\ldots
<i_{r}>i_{r+1}\neq \ldots \neq i_{n}$.
Since $(i_{1}, \ldots , i_{n})$ is non-crossing, hence there exists $1\leq j \leq r$ such that $i_j\neq i_k$
for every $k\neq j$. Then
\begin{eqnarray*}
\phi(a_{1}\ldots a_{n})&=&\phi(a_{1}\ldots a_{j-1}a_{j}^0a_{j+1}\ldots a_{n})
+\phi(a_{j})\phi(a_{1}\ldots a_{j-1}1_{i_j}a_{j+1}\ldots a_{n})\\
&=&
\phi(a_{j})\phi(a_{1}\ldots a_{j-1}a_{j+1}\ldots a_{n})
\end{eqnarray*}
by Lemma 3.5 since $a_{j}^0$ is a singleton of mean zero. To complete the proof, it is now enough to use
the induction argument since $(i_1, \ldots , i_{j-1}, i_{j+1}, \ldots , i_{n})\sim {\cal ONC}_{n-1}(m)$.
\hfill {$\blacksquare$}\\
\indent{\par}
{\sc Lemma 6.4.}
{\it Suppose that $i_{1}\neq i_{2}\neq \ldots \neq i_{n}$ and that
$(i_{1},\ldots , i_{n})\sim P\in {\cal ONC}_{n}\setminus {\cal ONC}_{n}(m)$.
Then the mixed moment $\phi(a_{1}\ldots a_{n})$ is a product of at least
$p+1$ marginal moments, where $p$ is the number of blocks of $P$.}\\
\indent{\par}
{\it Proof.}
Consider first the case when there exists $1\leq r \leq n$ such that $i_{1}\neq \ldots \neq i_{m}<\ldots <i_{r}>i_{r+1}$
and that for every $1\leq j \leq r$ there exists $r+2\leq k \leq n$ such that $i_j=i_k$. Then the second summand in
the expression
$$
\phi(a_{1}\ldots a_{n})=\phi(a_{1}^0a_{2}\ldots a_{n})
+\phi(a_{1})\phi(a_{2}\ldots a_{n})
$$
is a product of at least $p+1$ factors since the partition associated with $(i_{2}, \ldots , i_{n})$ has $p$ blocks.
Hence it suffices to consider the first summand
$$
\phi(a_{1}^0a_{2}\ldots a_{n})=\phi(a_{1}^0a_{2}^0a_{3}\ldots a_{n})
+\phi(a_{2})\phi(a_{1}^0a_{3}\ldots a_{n})
$$
and apply the same argument, using the property of $m$-monotone units. This reduces the proof to
the case of $\phi(a_{1}^0\ldots a_{r}^0a_{r+1}\ldots a_{n})$, which  vanishes by Lemma 3.4. The general case
can always be reduced to that considered above by means of Lemma 3.5 and
Lemma 6.1.
\hfill {$\blacksquare$}\\[10pt]
\myownsection
\begin{center}
{\sc 7. Central limit theorem}
\end{center}
We will formulate now the central limit theorem for the normalized sums of $m$-monotone random variables.
Then we shall give a recurrence relation for the sequence of their Cauchy transforms.\\
\indent{\par}
{\sc Theorem 7.1.}
{\it Let $\{X_i\}_{i=1}^{\infty}$ be identically distributed $m$-monotone random variables
with respect to $\phi$ of mean zero and variance one, i.e. $\phi(X_i)=0$ and $\phi(X_i^2)=1$. Then}
$$
\lim_{N\rightarrow \infty}\phi\left(\left[\frac{X_1+X_2+...+X_N}{\sqrt{N}}\right]^{n}\right)=
\begin{cases}
0 & \mbox{if $n=2k+1$} \\
\cfrac{\left|{\cal ONC}^2_{2k}(m)\right|}{k!} & \mbox{if $n=2k$}\\
\end{cases}.
$$
\indent{\par}
{\it Proof.}
Since mixed moments are invariant under order-preserving injections, we get
\begin{eqnarray*}
\phi((X_{1}+\ldots + X_{N})^{n})&=&
\sum_{1\leq k_1, \ldots , k_n \leq N}
\phi(X_{k_1}\ldots X_{k_n})\\
&=& \sum_{P\in {\cal OP}_{n}}{N \choose b(P)}m(P)
\end{eqnarray*}
where ${\cal OP}_{n}$ denotes all ordered partitions of the set $\{1, \ldots , n \}$,
$$
m(P)=\phi(X_{k_1}\ldots X_{k_n})\;\; {\rm for} \;\;(k_1, \ldots , k_{n})\sim P
$$
and $b(P)$ denotes the number of blocks of $P$.

Let us now analyze the contribution from every ordered partition $P$ for large $N$.
If there is a singleton in $P$, then by Lemma 6.1, $m(P)=0$. Thus suppose that $P$ has no singletons.
If $n=2k+1$, then $b(P)\leq k$ and after dividing ${N \choose b(P)}$ by
$N^{k+1/2}$ we get a contribution  of order smaller or equal to $N^{-1/2}$.
If $n=2k$, then either $b(P)<k$, when the contribution is of order smaller or equal to $N^{-1}$, or
$b(P)=k$, when we get $1/b(P)! m(P)$, but only for those ordered pair-paritions for which $m(P)$ factorizes
into $k$ factors. Lemmas 6.2-6.4 imply that if $P \in {\cal ONC}_{2k}^{2}(m)$, then $m(P)$ factorizes into
$b(P)$ factors, whereas if $P \notin {\cal OP}_{2k}^{2}\setminus {\cal ONC}_{2k}^{2}(m)$,
then $m(P)$ factorizes into more than $b(P)$ factors. This implies
the latter do not contribute to the limit. In turn, the contribution from the former is $1/k!$ since variance
is assumed to be equal to one. This completes the proof. \hfill {$\blacksquare$}\\
\indent{\par}
{\it Remark.}
If $m=1$, then using Theorem 7.1, we get in the CLT
\begin{eqnarray*}
\mathit{M}{}^{(1)}_{\,n}&=&
\begin{cases}
\ \ \ \,0 & \text{if $n=2k+1$}\\
\cfrac{1}{2^k}\left( \displaystyle{\genfrac{}{}{0pt}{}{2k}{k}} \right) & \text{if $n=2k$}
\end{cases}
\end{eqnarray*}
the moments of the arcsine law supported on $[-\sqrt{2},\sqrt{2}]$.\\
\indent{\par}
{\sc Theorem 7.2.}
{\it The Cauchy transforms $G^{(m)}(z)$ of the limit laws obtained in the CLT
for the monotone hierarchy satisfy the recurrence}
$$
G^{(1)}(z)=\frac{1}{\sqrt{z^2-2}}\ \ ,\ \ \ \ G^{(m)}(z)=\frac{1}{z-G^{(m-1)}(z)}\ \ \ {\it for}\ \ \ {\it m}\geqslant 2
$$
{\it for} $z\in\mathbb{C}_+=\{z:\mathfrak{Im}z>0\}$. \\
\indent{\par}
{\it Proof.}
That $G^{(1)}(z)$ is the Cauchy transform of the arcsine law is well-known (see, for instance [M2]).
Thus let $m\geqslant 2$.
Denote by $N_{2n}(m)$ the cardinality of ${\cal ONC}_{2n}^2(m)$.
We will express $N_{2n+2}(m)$ in terms of $N_{2k}(m)$ for $k\leqslant n$.
From the set $\{2,4,...,2n+2\}$ choose an even number which will form a pair with 1 and denote it $2k$.
Using the fact that partitions from ${\cal ONC}_{2n+2}^2(m)$ are non-crossing, we can separate
blocks which are inner w.r.t. $\{1,2k\}$ from other blocks (subsets of $\{2k+1,...,2n+2\}$).
The collection of the latter forms a partition of $\{2k+1,...,2n+2\}$
associated with ${\cal ONC}_{2n+2-2k}^2(m)$ and there are exactly
$\left(\genfrac{}{}{0pt}{}{n+1}{n+1-k}\right) N_{2n+2-2k}(m)$ such partitions.
On the other hand, blocks which are inner w.r.t. $\{1,2k\}$ form a partition associated with
${\cal ONC}_{2k-2}^2(m-1)$, and there are $\left(\genfrac{}{}{0pt}{}{k}{k-1}\right) N_{2k-2}(m-1)$
ways of choosing them. Thus
$$
N_{2n+2}(m)
=\sum_{k=1}^{n+1} k\left(\genfrac{}{}{0pt}{}{n+1}{k}\right)
N_{2k-2}(m-1)N_{2n-2k+2}(m) .
$$
Dividing both sides of the equation by $(n+1)!$ and denoting the moments of the central limit laws by $M_n^{(m)}$,
we obtain
$$
\frac{N_{2n+2}(m)}{(n+1)!}\ =\sum_{k=1}^{n+1}\ \frac{N_{2k-2}(m-1)}{(k-1)!}\ \frac{N_{2n-2k+2}(m)}{(n-k+1)!}
$$
$$
M_{2n+2}^{(m)}\ =\ \sum_{k=1}^{n+1}\ \ M_{2k-2}^{(m-1)}\ \ M_{2n-2k+2}^{(m)}\ .
$$
This identity allows us to derive a recurrence for the Cauchy transforms:
\begin{eqnarray*}
G^{(m)}(z)\ &=& \sum_{n=0}^{\infty}\ M_{\,2n}^{(m)}\ z^{-2n-1}\
=\ \frac{1}{z}+\ \sum_{n=0}^{\infty}\ M_{2n+2}^{\ (m)}\ z^{-2n-3}\\
&=&
\frac{1}{z}+\frac{1}{z}\ \sum_{n=0}^{\infty}\ \sum_{k=1}^{n+1}\ M_{\,\,2k-2}^{(m-1)}\ z^{-2k+1}\
M_{2n-2k+2}^{\ (m)}\ z^{-2n+2k-3}\\
&=&
\frac{1}{z}+\frac{1}{z} \bigg(\sum_{n=0}^{\infty}\ M_{\,2n}^{(m)}\ z^{-2n-1}\bigg)
\bigg(\sum_{n=0}^{\infty}\ M_{\,2n}^{(m-1)}\ z^{-2n-1}\bigg)\\
&=&
\frac{1}{z}+\frac{G^{(m)}(z)G^{(m-1)}(z)}{z}\ ,
\end{eqnarray*}
which leads to
$$
G^{(m)}(z)=\frac{1}{z-G^{(m-1)}(z)}\ ,
$$
and this completes the proof.\hfill{$\blacksquare$}\\
\indent{\par}
The table given below gives even moments of lowest orders obtained in
the CLT for the monotone hierarchy as well as those of the Wigner law which corresponds to
$m=\infty$. Let us observe that for given $m$, moments of orders smaller or equal to $2m$
agree with the moments of the Wigner law.

\begin{center}
\begin{tabular}{|c|c|c|c|c|c|}                                           \hline
     & $n=2$ &$n=4$        &$n=6$        &$n=8$         &  $n=10$       \\ \hline
$m=1$&  1    & 3/2         &       5/2   &       35/8   &       63/8    \\ \hline
$m=2$&  1    & 2           &       9/2   &       21/2   &       199/8   \\ \hline
$m=3$&  1    & 2           &5            &       27/2   &       75/2    \\ \hline
$m=4$&  1    & 2           &5            &14            &       83/2    \\ \hline
                                                                           \hline
$m=\infty$&1 & 2           &5            &14            &42             \\ \hline
\end{tabular}
\end{center}
\begin{center}
{\it Table 1.}
Low-order moments of the $m$-monotone central limit laws.
\end{center}

\indent{\par}
{\it Example 7.1.}
If $m=2$, then the Cauchy transform of the central limit law takes the form
$$
G^{(2)}(z)\ =\ \cfrac{1}{z-\cfrac{1}{z\sqrt{1-\cfrac{2}{z^2}}}}\ .
$$
In order to compute the absolutely continuous part of the
limit measure for $m=2$ we use the Stjelties inversion formula:
\begin{eqnarray*}
f^{(2)}(x)&=&-\frac{1}{\pi}\lim_{y\rightarrow 0^{\tiny +}}\mathfrak{Im}\ G^{(2)}\,(x\,+\,iy)\\
&=&
\begin{cases}
\cfrac{1}{\pi}\ \cfrac{\sqrt{2-x^2}}{1+x^2(2-x^2)} & \ \text{for}\ |x|\leqslant\sqrt{2}\\
0 & \ \text{for}\ |x|>\sqrt{2}\\
\end{cases}\ .
\end{eqnarray*}
The singular part of the limit measure is obtained
by computing the residua at the real poles of $G^{(2)}(z)$,
namely
$$
{\rm Res}_{\,\sqrt{\sqrt{2}+1}}\ G^{(2)}(z)\ =\ {\rm Res}_{\, - \sqrt{\sqrt{2}+1}}\ G^{(2)}(z)\ =\ \cfrac{2-\sqrt{2}}{4}\ .
$$
\noindent
Figure 3 below shows both the absolutely continuous and the singular parts of the central limit law
for $m=2$.

Using similar computations and Mathematica packet for $m=3$, we obtain a measure
with the absolutely continuous part
$$
f^{(3)}(x)=
\begin{cases}
\cfrac{1}{\pi}\ \cfrac{\sqrt{2-x^2}}{x^2+(x^2-1)^2(2-x^2)} & \ \text{for}\ |x|\leqslant\sqrt{2}\\
0 & \ \text{for}\ |x|>\sqrt{2}\\
\end{cases}\ .
$$
and the singular part consisting of two atoms of masses approximately equal to $0,099$ each,
concentrated at $\pm 1,685$, see Figure 4. An analogous picture for $m=4$ is given in Figure 5.\\

\unitlength=1mm
\special{em:linewidth 0.3pt}
\linethickness{0.3pt}
\begin{picture}(120.00,52.00)(-18,5)

\put(18.3,12.5){\includegraphics[]{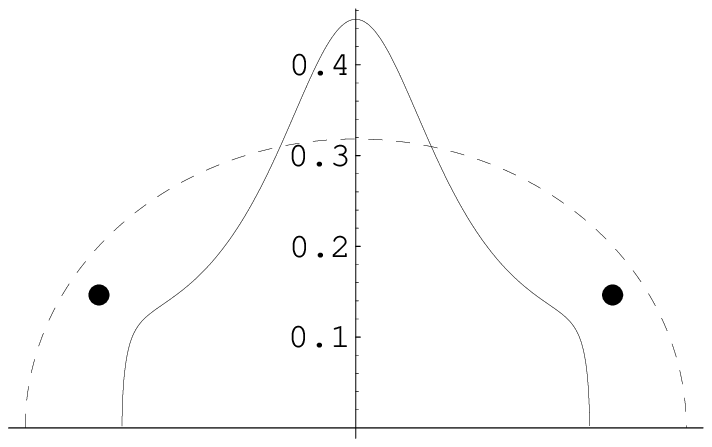}}

\put(75.5,10.9){\tiny{$\sqrt{2}$}}
\put(86.5,10.9){\tiny{$2$}}

\put(25.5,10.9){\tiny{$-\sqrt{2}$}}
\put(19.3,10.9){\tiny{$2$}}
\put(15,5){{\it Figure 3.} The central limit law for $m=2$.}
\end{picture}\\

\unitlength=1mm
\special{em:linewidth 0.3pt}
\linethickness{0.3pt}
\begin{picture}(120.00,50.00)(-18,7)
\put(18.3,12.5){\includegraphics[]{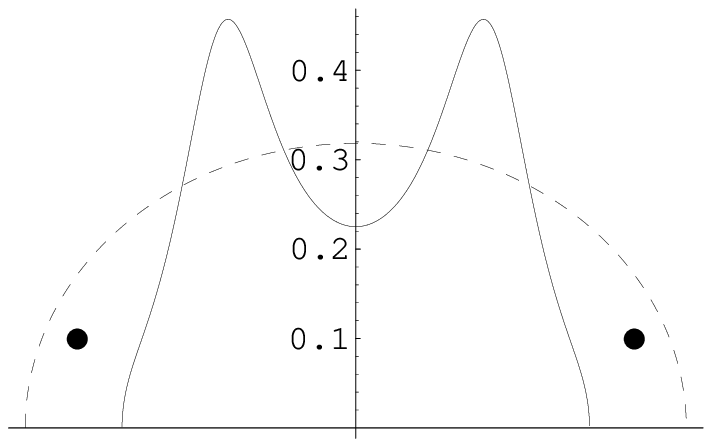}}

\put(75.5,10.9){\tiny{$\sqrt{2}$}}
\put(86.5,10.9){\tiny{$2$}}

\put(25.5,10.9){\tiny{$-\sqrt{2}$}}
\put(19.3,10.9){\tiny{$2$}}
\put(15,4){{\it Figure 4.} The central limit law for $m=3$.}
\end{picture}
\\[10pt]

\unitlength=1mm
\special{em:linewidth 0.3pt}
\linethickness{0.3pt}
\begin{picture}(120.00,50.00)(-18,7)
\put(18.3,12.5){\includegraphics[]{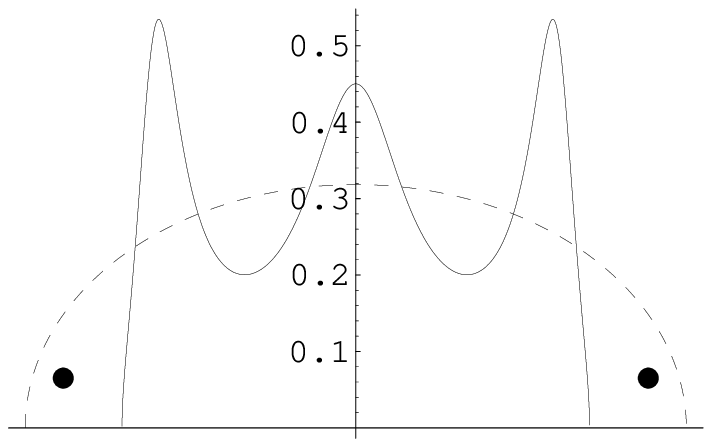}}

\put(75.5,10.9){\tiny{$\sqrt{2}$}}
\put(86.5,10.9){\tiny{$2$}}

\put(25.5,10.9){\tiny{$-\sqrt{2}$}}
\put(19.3,10.9){\tiny{$2$}}
\put(15,4){{\it Figure 5.} The central limit law for $m=4$.}
\end{picture}
\\[10pt]

\myownsection
\begin{center}
{\sc 8. Poisson's limit theorem.}
\end{center}
In this Section we prove the Poisson's limit theorem. \\
\indent{\par}
{\sc Theorem 8.1.}
{\it Let $m\in {\mathbb N}$ and suppose that for every $N\in {\mathbb N}$
the variables $X_{1,N}, \ldots , X_{N,N}$ are identically distributed
and $m$-monotone with respect to $\phi_{N}$.
If $N\phi_{N}(X_{i,N}^{k})\rightarrow \lambda^{k}$ for every natural $k$, where $\lambda>0$, then}
$$
\lim_{N\rightarrow \infty}\phi_{N}((X_{1,N}+ \ldots + X_{N,N})^{n})=
\sum_{q=1}^{n}\frac{\lambda^{q}}{q!}|{\cal ONC}_{n}(q,m)|
$$
{\it where ${\cal ONC}_{n}(q,m)$ denotes the subset of ${\cal ONC}_{n}(m)$ consisting
of partitions which have $q$ blocks.}\\
\indent{\par}
{\it Proof.}
Using invariance of mixed moments under order-preserving injections, we get
\begin{eqnarray*}
\phi_{N}((X_{1,N}+\ldots + X_{N,N})^{n})&=&
\sum_{1\leq k_1, \ldots , k_n \leq N}\phi_{N}(X_{k_1,N}\ldots X_{k_n,N})\\
&=& \sum_{P\in {\cal OP}_{n}}  {N \choose b(P)}m_{N}(P)
\end{eqnarray*}
where ${\cal OP}_{n}$ denotes all ordered partitions of the set $\{1, \ldots , n \}$,
$$
m_{N}(P)=\phi_{N}(X_{k_1,N}\ldots X_{k_n,N})\;\; {\rm whenever} \;\;(k_1, \ldots , k_{n})\sim P
$$
and $b(P)$ denotes the number of blocks of $P$.

Lemmas 6.2-6.4 imply that if $P \in {\cal ONC}_{n}(m)$, then $m_N(P)$ factorizes into $b(P)$ factors,
whereas if $P \notin {\cal ONC}_{n}(m)$, then $m(P)$ factorizes into more than $b(P)$ factors. This implies
that the latter do not contribute to the limit. In turn, the contribution from $P\in {\cal ONC}_{n}(m)$,
depending on $b(P)$, gives the desired formula.\hfill {$\blacksquare$}\\[10pt]
\indent{\par}
Again, as in the case of the central limit theorem, the recurrence for Cauchy transforms
is the same as in the $m$-free case [L-F] and only the intitial condition is different.
We use the notations:
$$
M_{n}^{(m)}(\lambda,q)=\frac{\lambda^{q}}{q!}|{\cal ONC}_{n}(q,m)|\ \ \ ,\ \ \ \ \ H^{(m)}(\lambda,z)=\sum_{n,q=0}^{\infty} M_{n}^{(m)}(\lambda,q)\, z^{-n-1}
$$
with $M_{n}^{(m)}(\lambda,0)=\delta_{n,0}$.\\
\indent{\par}
{\sc Lemma 8.2.}
{\it The hierarchy of generating functions $\{H^{(m)}\}_{m=0}^{\infty}$ satisfies the recurrence relation}
$$
H^{(m)}(\lambda,z)=\frac{1-H^{(m-1)}(\lambda,z)}{z-zH^{(m-1)}(\lambda,z)-\lambda}
$$
{\it for $m=2,3,...$ and $H^{(1)}(\lambda,z)$ is the product log function described in} [M2].\\
\indent{\par}
{\it Proof.} The case $m=1$ follows from [M2]. Suppose that $m\geqslant 2$.
To get an ordered non-crossing partition of $\{1,...,n\}$ we pick the elements
that will be put in the same block as the first element; denote this block by
$\{1,p_{1}, \ldots , p_{r-1}\}$ as in the diagram given below, with
$Q_{1},Q_{2}, \ldots , Q_{r}$ denoting subpartitions of the intervals
$(1,p_{1}-1)$, $(p_{1}+1, \ldots , p_{2}-1), \ldots  (p_{r-1},n]$, respectively.\\
\unitlength=1mm
\special{em:linewidth 0.4pt}
\linethickness{0.4pt}
\begin{picture}(120.00,42.00)(-35.00,5.00)

\put(0.00,21.00){\line(0,1){12.00}}
\put(15.00,21.00){\line(0,1){12.00}}
\put(30.00,21.00){\line(0,1){12.00}}
\put(0.00,33.00){\line(1,0){40.00}}
\put(60.00,21.00){\line(0,1){12.00}}
\put(75.00,21.00){\line(0,1){12.00}}
\put(55.00,33.00){\line(1,0){20.00}}

\put(0.00,21.00){\circle*{1.00}}
\put(15.00,21.00){\circle*{1.00}}
\put(30.00,21.00){\circle*{1.00}}
\put(60.00,21.00){\circle*{1.00}}
\put(75.00,21.00){\circle*{1.00}}

\put(45.00,33.00){...}
\put(45.00,21.00){...}
\put(-1.00,17.00){\footnotesize $1$}
\put(14.00,17.00){\footnotesize $p_{1}$}
\put(29.00,17.00){\footnotesize $p_{2}$}
\put(57.00,17.00){\footnotesize $p_{r-2}$}
\put(72.00,17.00){\footnotesize $p_{r-1}$}

\put(7.00,21.00) {\footnotesize $Q_{1}$}
\put(22.00,21.00) {\footnotesize $Q_{2}$}
\put(62.00,21.00) {\footnotesize $Q_{r-1}$}
\put(85.00,21.00) {\footnotesize $Q_{r}$}

\end{picture}
\\
As compared with the $m$-free case we need to take into account the number
of possible block colorings. Denoting $N_{n}(q,m)=|{\cal ONC}_{n}(q,m)|$, we get
\begin{eqnarray*}
N_{n}(q,m)&=&
\sum_{r=1}^{n}\ \sum_{k_1+...k_{r}=n-1}\ \sum_{q_1+...+q_r=q-1}\
\frac{q!}{q_1!q_2!...q_r!}\\
&\times&
N_{k_1-1}(q_1,m-1)... N_{k_{r-1}-1}(q_{r-1},m-1)N_{k_r}(q_r,m)
\end{eqnarray*}
for $n,k_1,...,k_{r-1}\geqslant1$\,; \ $k_r,b_1,...,b_r\geqslant0$,
where
$$
k_{1}=p_{1}-1,k_{2}=p_{2}-p_{1}, \ldots , k_{r-1}=p_{r-1}-p_{r-2}
$$
since there are $q!/(q_{1}!\ldots q_{r}!)$ ways to choose {\it sets} of colors
for subparitions $Q_{1},\ldots , Q_{r}$ from among $q=q_{1}+\ldots q_{r}$ colors and
then `admissible' colorings of these sets give numbers
$N_{k_1-1}(q_1,m-1)$, $\ldots $, $N_{k_r}(q_r,m)$.
Now, multiplying both sides by $\lambda^q/q!$, we get
\begin{eqnarray*}
M_{n}^{(m)}(\lambda,q)&=&\lambda \sum_{r=1}^{n}\ \sum_{k_1+...+k_r=n-1}\
\sum_{q_1+...+q_r=q-1}\\
&&M_{k_1-1}^{(m-1)}(\lambda,q_1)
\times...M_{k_{r-1}-1}^{(m-1)}(\lambda,q_{r-1})\, M_{k_1}^{(m)}(\lambda,q_r)
\end{eqnarray*}
Then we have
\begin{eqnarray*}
H^{(m)}(\lambda,z)&=&\sum_{n,q=0}^{\infty} M_{n}^{(m)}(\lambda,q)\, z^{-n-1}\\
&=&\frac{1}{z}+\frac{\lambda}{z}\sum_{n,q=1}^{\infty}\,
\sum_{r=1}^{n}\ \sum_{k_1+...+k_r=n-1}\ \sum_{q_1+...+q_r=q-1}\
M_{k_1-1}^{(m-1)}(\lambda,q_1)\,z^{-k_1}\\
& &\times...M_{k_{r-1}-1}^{(m-1)}(\lambda,q_{r-1})\,z^{-k_{r-1}}\
M_{k_1}^{(m)}(\lambda,q_r)\,z^{-k_r-1}\\
&=&\frac{1}{z}+\frac{\lambda}{z}\sum_{r=1}^{\infty}
\Big( \sum_{\beta,\nu=0}^{\infty} M_{\nu}^{(m-1)}(\lambda,\beta)\,z^{-\nu -1}\Big)^{r-1}
\sum_{\mu,\alpha=0}^{\infty} M_{\mu}^{(m)}(\lambda,\alpha)\,z^{-\mu -1}\\
&=&\frac{1}{z}+\frac{\lambda H^{(m)}(\lambda,z)}{z\big(1-H^{(m-1)}(\lambda,z)\big)}\ ,\\
\end{eqnarray*}

\noindent
and therefore
$$
H^{(m)}(\lambda,z)=\frac{1-H^{(m-1)}(\lambda,z)}{z-zH^{(m-1)}(\lambda,z)-\lambda}
$$
which completes the proof.\hfill {$\blacksquare$}\\[10pt]

\myownsection
\begin{center}
{\sc 9. Invariance principle}
\end{center}
Let $(X_{i})_{i\in {\mathbb N}}\in {\cal A}$ be a sequence of
identically distributed random variables of mean zero and variance one, which are $m$-monotone with respect
to a state $\phi$ on ${\cal A}$. Below we shall investigate the
asymptotic behavior of normalized sums
\begin{equation}
S_{N}(f)=\frac{1}{\sqrt{N}}\sum_{i=[sN]+1}^{[tN]}X_{i}
\end{equation}
for $f=\chi{(s,t]}\in \Theta$ (see (5.1)) as  $N\rightarrow \infty$.

It is convenient to give a formula for mixed moments of sums indexed by functions
which have identical or disjoint supports. \\
\indent{\par}
{\sc Theorem 9.1.}
{\it Suppose $f_{1}, f_{2},\ldots , f_{n}\in \Theta$ have identical or disjoint supports
and let $\sigma(f) =\{\sigma_{1}, \sigma_{2}, \ldots , \sigma_{p}\}$ be the partition defined
by the tuple $(f_{1},f_{2}, \ldots , f_{n})$, i.e. each block $\sigma_{k}$ consists of
all those indices $j$, for which $f_{j}$ coincide. Then}
\begin{equation}
\lim_{N\rightarrow \infty}
\phi(S_{N}(f_{1})S_{N}(f_{2})\ldots S_{N}(f_{n}))
=
\frac{1}{b_1!b_2!\ldots b_p!}
\sum_{
\stackrel{P\in {\cal ONC}_{n}^{2}(m)}
{\scriptscriptstyle P\sim (f_1,f_2, \ldots ,f_n)}
}
\prod_{\{\alpha , \beta\}\in P }
\langle f_{\alpha}, f_{\beta} \rangle
\end{equation}
{\it where $b_{k}=|\sigma_{k}|/2$ for $k=1,\ldots , p$, with the understanding that ${\cal ONC}_{n}^{2}(m)=\emptyset$
for $n$ odd.}\\
\indent{\par}
{\it Proof.}
It is clear that if $n$ is odd, we get zero contribution - we understand here that in that case
${\cal ONC}_{n}^{2}(m)=\emptyset$. Thus, let $n=2k$ and let $f_{k}=\chi_{(s_{k},t_{k}]}$ for
$k=1, \ldots , n$. We have
\begin{eqnarray*}
\phi(S_{N}(f_{1})S_{N}(f_{2})\ldots S_{N}(f_{n}))
&=&
\frac{1}{N^{n/2}}
\sum_{i_1=[s_1N]+1}^{[t_1N]}
\ldots
\sum_{i=[s_nN]+1}^{[t_nN]}
\phi(X_{i_1}\ldots X_{i_n})\\
&=&
\frac{1}{N^{n/2}}
\sum_{P\in {\cal ONC}_{n}^{2}(m)}A_{P}(f_{1},f_{2}, \ldots , f_{n};N)+{\cal O}(1/N)
\end{eqnarray*}
where $A_{P}(f_{1},f_{2}, \ldots , f_{n};N)$ denotes the number of tuples $(i_{1},i_{2}, \ldots , i_{n})$
such that
$$
[s_jN]<i_{j}\leq [t_jN] \;\;{\rm for}\;\; j=1, \ldots , n
$$
and $(i_{1},i_{2}, \ldots , i_{n})\sim P$.
Now, Lemma 6.1 together with the mean zero assumption result in the fact that only contributions from
pair partitions are of order $1$ - other partitions must have fewer than $k$ blocks and thus their contribution is at
most of order $1/N$ as in the CLT. Moreover, moments associated with
ordered pair partitions which are not in ${\cal ONC}_{n}^{2}(m)$ factorize into
more than $k$ factors by Lemmas 6.2-6.4 and thus vanish by the mean zero assumption.
Finally, if $P\in {\cal ONC}_{n}^{2}(m)$ and $P\sim (f_{1},f_{2}, \ldots , f_{n})$, then
$$
A_{P}(f_{1}, f_{2},  \ldots ,f_{n};N)
=
{N_{1}\choose b_{1}}
{N_{2}\choose b_{2}}
\ldots
{N_{p}\choose b_{p}}
$$
where $N_{j}=[t^{(j)}N]-[s^{(j)}N]$ and $t^{(j)}, s^{(j)}$ denote the end-points of $f_r$ for any
$r\in \sigma_{j}$, where $j=1, \ldots , p$
(thus $N_{j}$ counts the number of possible values of those
indices in the tuple $(i_1,i_2, \ldots ,i_n)$ which are associated with block $\sigma_j$).
For such $P$ we get
\begin{eqnarray*}
\frac{A_{P}(f_{1},f_{2}, \ldots , f_{n})}
{N^{k}}
&=&
\prod_{j=1}^{p}\frac{N_{j}(N_{j}-1)\ldots (N_{j}-b_{j}+1)}{b_{j}!N^{b_{j}}}\\
&\rightarrow&
\prod_{j=1}^{p}\frac{(t^{(j)}-s^{(j)})^{b_{j}}}{b_{j}!}\\
&=&
\frac{1}{b_{1}!b_{2}!\ldots b_{p}!}\prod_{\{\alpha , \beta \}\in P}
\langle f_{\alpha}, f_{\beta} \rangle
\end{eqnarray*}
whereas if $P\not \sim (f_{1},f_{2},\ldots , f_{n})$, then $A_{P}(f_{1},f_{2}, \ldots , f_{n};N)=0$,
which finishes the proof.
\hfill {$\blacksquare$}\\[10pt]

\myownsection
\begin{center}
{\sc 10. $m$-monotone Gaussian operators}
\end{center}
In this section we introduce the $m$-monotone Fock space as an $m$-th order iteration
of the free Fock space, on which we give a realization of the limit moments in the invariance principle
for $m$-monotone independence.\\

The $m$-{\it free Fock space} over a Hilbert space ${\cal H}$ [FL] is defined by the direct sum
$$
{\cal F}_{(m)}({\cal H})={\mathbb C}\Omega \oplus \bigoplus_{n=1}^{m}{\cal H}^{\otimes n}
$$
with the canonical scalar product inherited from the free Fock space ${\cal F}({\cal H})$, where
$m\in {\mathbb N}$ and $\Omega$ is a unit vector called vacuum.

Let ${\cal H}=L^{2}({\mathbb R}_{+})$, then ${\cal H}^{\otimes n}\cong L^{2}({\mathbb R}_{+}^{n})$
and thus another direct sum
\begin{equation}
{\cal F}^{(m)}({\cal H})={\cal F}_{(m)}({\cal H}) \oplus \bigoplus_{n=m+1}^{\infty}L^{2}(\Delta_{n}^{(m)})
\end{equation}
can be viewed as a subspace of ${\cal F}({\mathbb R}_{+})$ with `order starting from the $m$-th level', where
$$
\Delta_{n}^{(m)}=\{(x_{1},x_{2}, \ldots , x_{n}):x_{1}>\ldots > x_{n-m+1}>0\}
$$
for every $n\in {\mathbb N}$. We equip ${\cal F}^{(m)}({\cal H})$
with the canonical scalar product inherited from ${\cal F}({\cal H})$ and call it
the $m$-{\it monotone Fock space}. In particular, ${\cal F}^{(1)}({\cal H})$ is the
monotone Fock space over ${\cal H}$ [M2,Lu].

If we denote by ${\cal P}^{(m)}$ the canonical projection from
${\cal F}({\cal H})$ onto ${\cal F}^{(m)}({\cal H})$, then ${\cal F}^{(m)}({\cal H})$
is spanned by $\Omega$ and vectors of the form
$$
f_{1}\otimes _{m}\ldots \otimes _{m} f_{n}
:={\cal P}^{(m)}(f_{1} \otimes \ldots \otimes f_{n})
$$
where $f_{1}, \ldots , f_{n}\in {\cal H}$ and $n\in {\mathbb N}$.

Define $m$-{\it monotone creation operators} on ${\cal F}^{(m)}({\cal H})$ as bounded extensions of
\begin{eqnarray*}
a^{(m)}(f)\Omega&=&f\\
a^{(m)}(f)(f_{1}\otimes _{m}\ldots  \otimes _{m} f_{n})&=&
f\otimes _{m}f_{1} \otimes _{m}\ldots \otimes _{m} f_{n}
\end{eqnarray*}
and the $m$-{\it monotone annihilation operators} $a^{(m)*}(f)$ as their adjoints.

It is easy to see that $a^{(m)*}(f)\Omega =0$ and that
$$
a^{(m)*}(f) (f_{1}\otimes _{m}\ldots \otimes _{m} f_{n})
=\left\{
\begin{array}{ll}
\langle f_{1},f\rangle f_{2}\otimes _{m}\ldots \otimes _{m} f_{n} & 1\leq n \leq m\\
M_{\psi}f_{2}\otimes _{m}\ldots \otimes _{m}f_{n} & n>m
\end{array}
\right.
$$
where $(M_{\psi}g)(x)=\psi(x)g(x)$ is the multiplication operator by function $\psi$
and
\begin{equation}
\psi(x)=\langle f_{1},f \rangle _{x}= \int_{y>x}f_{1}(y)\bar{f}(y)dy
\end{equation}
(dependence of $\psi$ on $f$ and $f_{1}$ is supressed).
Note that if ${\rm supp}f_{1}\cap {\rm supp}f > {\rm supp}f_{2}$, then
we have $\psi(x)=\langle f_{1},f\rangle $ for every $x\in {\rm supp}f_{2}$.

Now let us define vacuum states on all bounded operators $T$ on the Fock spaces
${\cal F}^{(m)}({\cal H})$:
$$
\varphi (T)=\langle T\Omega , \Omega \rangle
$$
(by abuse of notation we use the same notation for all $m$), and position operators
$$
\omega^{(m)}(f)=a^{(m)}(f)+a^{(m)*}(f)
$$
for $f\in \Theta$ given by (5.1), which may be called $m$-{\it monotone Gaussian operators}.

Our goal is now to express the limit moments obtained
in the invariance principle in terms of the vacuum expectations of
the $m$-monotone Gausssian operators.
For clarity of exposition, we first consider the monotone case $m=1$.
A generalization to all natural $m$ is rather straightforward.
The combinatorial formulas for mixed moments even in the case of monotone Gaussians
(with arcsine distribution) are new.

To every $\pi\in {\cal NC}_{2k}^{2}$ and $(f_{1},f_{2} \ldots , f_{2k})$  we associate the operator
\begin{equation}
a_{\pi}(f_{1},f_{2}, \ldots , f_{2k})=
a^{\epsilon_{1}}(f_{1})a^{\epsilon_{2}}(f_{2})\ldots a^{\epsilon_{2k}}(f_{2k})
\end{equation}
where $\epsilon_{p}=*$ and $\epsilon_{q}=1$ whenever $\{p,q\}$ is a block
of $\pi$. We also use notations (5.2) and (5.3).\\
\indent{\par}
{\sc Lemma 10.1.}
{\it Let $\pi=\{\pi_{1},\pi_{2}, \ldots , \pi_{k}\}\in {\cal NC}_{2k}^{2}$ and
$\pi\sim (f_{1}, f_{2}, \ldots , f_{2k})$, where $f_{1}, f_{2}, \ldots , f_{n}\in {\Theta}$
have identical or disjoint supports. Then}
\begin{equation}
\varphi (a_{\pi}(f_{1},f_{2}, \ldots ,f_{2k}))
=
\psi_{\pi_{1}}(s^{(1)})\psi_{\pi_{2}}(s^{(2)})\ldots \psi_{\pi_{k}}(s^{(k)})
\end{equation}
{\it where $\psi_{\pi_{i}}(x)=({\rm Inn}(\pi_{i})+1)^{-1}(t^{(i)}-x)$ is
the linear function associated with block $\pi_{i}$ for $x\in [s^{(i)},t^{(i)}]$,
$i=1, \ldots , k$.}\\
\indent{\par}
{\it Proof.}
Clearly, if $k=1$ and $\pi=\{\{1,2\}\}$, then
$$
a_{\pi}\Omega=a^{*}(f_{1})a(f_{1})\Omega = (t^{(1)}-s^{(1)})\Omega =\psi_{\pi_{1}}(s^{(1)})\Omega,
$$
which gives the formula for $k=1$.

{\bf Claim:}
Abbreviating notation (10.3) to $a_{\pi}$, we get
$$
a_{\pi}(g_{1}\otimes_{1} \ldots \otimes_{1} g_{n})=
\prod_{
\stackrel{1 \leq i\leq k}
{\scriptscriptstyle
{\rm supp}f_{i}\neq {\rm supp g}_{1}}}
\psi_{\pi_{i}}(s^{(i)})
(
\prod_{
\stackrel{1 \leq j\leq k}
{\scriptscriptstyle
{\rm supp}f_{j}={\rm supp g}_{1}}}
M_{\psi_{j}}g_{1}
)
\otimes_{1}\ldots \otimes_{1} g_{n}
$$
for $g_{1}, \ldots , g_{n}\in {\cal H}$
with ${\rm supp}g_{1}={\rm supp}f_{2k}$ or ${\rm supp}g_{1}<{\rm supp}f_{2k}$,
where $n\in {\mathbb N}$.

Suppose this formula holds for $\pi\in {\cal NC}_{2r}^{2}$ for $1\leq r \leq k-1$.
If $\pi=\pi'\cup \pi''$, where
$\pi'$ is a non-crossing partition of $\{1, \ldots , 2r\}$
and $\pi''$ is a non-crossing partition of $\{2r+1, \ldots , 2k\}$,
then it is obvious that the claim holds for $\pi$.
Therefore, suppose that $\pi=\pi'\cup \{\{1,2k\}\}$.
Then, using the inductive assumption (and abbreviated self-explanatory notation $a_{\pi'}$),
we get
\begin{eqnarray*}
a_{\pi}(g_{1}\otimes_{1}\ldots \otimes_{1}g_{n})
&=& a^{*}(f_{1})a_{\pi'}(f_{2k}\otimes_{1}g_{1}\otimes_{1}\ldots \otimes_{1} g_{n})\\
&=&
\prod_{
\stackrel{2 \leq j\leq k}
{\scriptscriptstyle f_{j}\neq f_{2k}}}
\psi_{\pi_{i}}(s^{(i)})
a^{*}(f_{1})
(
\prod_{
\stackrel{2 \leq j\leq k}
{\scriptscriptstyle f_{j}=f_{2k}}}
M_{\psi_{j}}f_{2k}\otimes_{1}g_{1} \otimes_{1}\ldots \otimes_{1}g_{n})\\
&=&
\prod_{
\stackrel{2 \leq j\leq k}
{\scriptscriptstyle f_{i}\neq f_{2k}}}
\psi_{\pi_{i}}(s^{(i)})
(M_{\psi}g_{1} \otimes_{1}\ldots \otimes_{1} g_{n}),
\end{eqnarray*}
where
$$
\psi (x)=
\left\{
\begin{array}{ll}
\psi_{\pi_1}(x) & {\rm if}\;{\rm supp}g={\rm supp}f_{2k}\\
\psi_{\pi_1}(s^{(1)}) & {\rm if} \;{\rm supp}g<{\rm supp}f_{2k}
\end{array}
\right.
$$
follows from
$$
\psi(x)=\int_{x}^{t^{(1)}}
\prod_{\stackrel
{2\leq j\leq k}
{\scriptscriptstyle  f_{j}= f_{2k}}
}
M_{\psi_{j}}1dx
=
\int_{x}^{t^{(1)}}
(t^{(1)}-y)^{b_{1}-1}dy
=
\frac{1}{b_{1}}(t^{(1)}-x)^{b_{1}}
$$
where $b_{1}={\rm Inn}(\pi_{1}) +1$, whenever ${\rm supp}g={\rm supp}f_{2k}$, whereas if
${\rm supp}g\neq {\rm supp}f_{2k}$, then $\psi(x)=\psi_{\pi_{1}}(s^{(1)})$.
This proves the claim, from which (10.8) follows easily.
\hfill {$\blacksquare$}\\
\indent{\par}
{\it Example 10.1.}
Let us evaluate the moment $\langle a_{\pi}\Omega, \Omega \rangle $, where
$\pi$ is given by the diagram in Figure 2 (see Section 5).
In view of Lemma 10.1, we get
$$
\langle a_{\pi}\Omega, \Omega \rangle =\frac{1}{2}t^{2}(t'-t)^{2}
$$
since
$\phi_{1}(x)=1/2(t-x)$, $\phi_{2}(x)=t'-x$, $\phi_{3}(x)=t-x$, $\phi_{4}(x)=t'-x$
(cf. Example 5.3).\\
\indent{\par}
For comparison, note that in the case of vacuum expectations of
free creation and annihilation operators
there is a multiplicative formula similar to (10.4), except that on the RHS we get
$(t^{(1)}-s^{(1)})\ldots (t^{(k)}-s^{(k)})$. Thus, the difference is that
in the monotone case the product encodes additional information about inner blocks.

Moreover, as the lemma below demonstrates,
the number of admissible colorings of $\pi$ can also be expressed
in terms of the same information about inner blocks.\\
\indent{\par}
{\sc Lemma 10.2.}
{\it Let $\pi=\{\pi_{1},\pi_{2}, \ldots , \pi_{k}\}\in {\cal NC}_{2k}^{2}$ and
$\pi\sim (f_{1}, f_{2}, \ldots , f_{2k})$, where
$f_{1}, f_{2}, \ldots , f_{n}\in {\Theta}$ have identical or disjoint supports.
Then}
\begin{equation}
\frac{c_{\pi}(f_{1}, f_{2}, \ldots , f_{2k})}
{b_{1}!b_{2}!\ldots b_{p}! }
=
\prod_{i=1}^{k}({\rm Inn}(\pi_{i})+1)^{-1}
\end{equation}
{\it with the same notation as in Theorem 9.1.}\\
\indent{\par}
{\it Proof.}
Clearly, the formula is true for $k=1$ and $\pi=\{1,2\}$. Suppose the formula holds for
non-crossing pair partitions of sets consisting of $2k-2$ elements and let
$\pi \in {\cal NC}_{2k}^{2}$.

{\it Case 1.}
Suppose $\pi=\pi'\cup \pi''$, where $\pi'=\{\pi_{1}', \ldots , \pi_{r}'\}$
is a non-crossing pair partition of $\{1, \ldots , 2r\}$
and $\pi''=\{\pi''_{1}, \ldots , \pi''_{k-r}\}$
is a non-crossing pair-partition of $\{2r+1, \ldots , 2k\}$.
Let $a_{i}$ denote the number of non-empty blocks of the partition
$\tau^{(i)}=\{\pi'_{1}\cap \sigma_{i}, \ldots , \pi'_{r}\cap \sigma_{i}\}$,
where $i=1, \ldots , p$.
Since $\pi\sim (f_{1}, \ldots , f_{2k})$, there are ${b_{i}\choose a_{i}}$ ways to choose colors
for blocks of $\tau^{(i)}$. Since colors are chosen independently for different supports, we get
\begin{eqnarray*}
\frac{c_{\pi}(f_{1}, \ldots , f_{2k})}
{b_{1}!\ldots b_{p}! }
&=&
{b_{1}\choose a_{1}}
\ldots
{b_{p}\choose a_{p}}
\frac
{c_{\pi'}(f_{1}, \ldots , f_{2r})c_{\pi''}(f_{2r+1}, \ldots , f_{2k})}
{b_{1}!\ldots b_{p}!}\\
&=&
\frac{c_{\pi'}(f_{1}, \ldots , f_{2r})}
{a_{1}!\ldots a_{p}! }
\frac{c_{\pi''}(f_{2r+1}, \ldots , f_{2k})}
{(b_{1}-a_{1})!\ldots (b_{p}-a_{p})! }\\
&=&
\prod_{i=1}^{r}({\rm Inn}(\pi'_{i})+1)^{-1}\prod_{j=1}^{k-r}({\rm Inn}(\pi''_{j})+1)^{-1}\\
&=&
\prod_{i=1}^{k}({\rm Inn}(\pi_{i})+1)^{-1}
\end{eqnarray*}
which gives (10.9) for the partition $\pi$.

{\it Case 2.}
Suppose $\pi=\pi'\cup \{\{1,2k\}\}$, where $\pi'=\{\pi_{2}, \ldots \pi_{k}\}$ is a non-crossing pair-partition
of the set $\{2, \ldots , 2k-1\}$. Thus all blocks of $\pi'$ are inner w.r.t. the block $\pi_{1}=\{1,2k\}$
and this implies that $\pi_{1}$ must be colored by $1$. Therefore, there are exactly the same numbers
of ways of coloring $\pi$ as well as $\pi'$.
Choosing $b_{1}$ to be the number of blocks of $\pi$ with the same support as $f_{1}=f_{2k}$, we get
\begin{eqnarray*}
c_{\pi}(f_{1}, \ldots , f_{2k})
&=&
c_{\pi'}(f_{2}, \ldots , f_{2k-1})\\
&=&
(b_{1}-1)!b_{2}!\ldots b_{p}!
\prod_{i=2}^{k}({\rm Inn}(\pi_{i})+1)^{-1}\\
&=&
b_{1}!b_{2}!\ldots b_{p}!
\prod_{i=1}^{k}({\rm Inn}(\pi_{i})+1)^{-1}
\end{eqnarray*}
since $b_{1}={\rm Inn}(\pi_{1})+1$, which gives (10.9) for the partition $\pi$. This completes the proof.
\hfill {$\blacksquare$}\\
\indent{\par}
{\sc Theorem 10.3.}
{\it Suppose $f_{1}, f_{2}, \ldots , f_{n}\in {\Theta}$ have identical or disjoint supports.
Then}
$$
\varphi(\omega^{(1)}(f_{1})\omega^{(1)}(f_{2})\ldots \omega^{(1)}(f_{n}))
=
\frac{1}{b_1!b_2!\ldots b_p!}
\sum_{
\stackrel{P\in {\cal MNC}_{n}^{2}}
{\scriptscriptstyle P\sim (f_1,f_2, \ldots ,f_n)}
}
\prod_{\{\alpha , \beta\}\in P }
\langle f_{\alpha}, f_{\beta} \rangle
$$
{\it with the same notation as in Theorem 9.2.} \\
\indent
{\it Proof.}
If $n=2k$, we have
\begin{eqnarray*}
\varphi(\omega^{(1)}(f_{1})\ldots \omega^{(1)}(f_{2k}))
&=&
\sum_{
\stackrel{\pi\in {\cal NC}_{2k}^{2}}
{\scriptscriptstyle \pi \sim (f_1, \ldots ,f_n)}
}
\langle a_{\pi}(f_{1}, \ldots , f_{2k})\Omega , \Omega \rangle\\
&=&
\sum_{
\stackrel{\pi\in {\cal NC}_{2k}^{2}}
{\scriptscriptstyle \pi \sim (f_1, \ldots ,f_n)}
}
\psi_{\pi_{1}}(s^{(1)})\ldots \psi_{\pi_{k}}(s^{(k)})\\
&=&
\sum_{
\stackrel{\pi\in {\cal NC}_{2k}^{2}}
{\scriptscriptstyle \pi \sim (f_1, \ldots ,f_n)}
}
\prod_{i=1}^{k}\frac{t^{(i)}-s^{(i)}}{{\rm Inn}(\pi_{i})+1}\\
&=&
\frac{1}{b_{1}!\ldots b_{p}!}
\sum_{
\stackrel{\pi\in {\cal NC}_{2k}^{2}}
{\scriptscriptstyle \pi \sim (f_1, \ldots ,f_n)}
}
c_{\pi}(f_{1}, \ldots , f_{2k})
\prod_{\{\alpha, \beta\}\in \pi}
\langle f_{\alpha}, f_{\beta}\rangle\\
&=&
\frac{1}{b_{1}!\ldots b_{p}!}
\sum_{
\stackrel{P  \in {\cal MNC}_{2k}^{2}}
{\scriptscriptstyle P \sim (f_1, \ldots ,f_n)}
}
\prod_{\{\alpha, \beta\}\in P}
\langle f_{\alpha}, f_{\beta}\rangle
\end{eqnarray*}
where we used Lemmas 10.2-10.3. It is clear that if $n$ is odd, then the expectation vanishes.
This completes the proof.
\hfill {$\blacksquare$} \\
\indent{\par}
In particular, the combinatorics of inner blocks allows us to express
the moments of the arcsine law as a sum
$M_{2k}^{(1)}=\sum_{\pi\in{\cal NC}_{2k}^{2}} \varphi(a_{\pi})$,
where the vacuum expectations $\varphi(a_{\pi})$ are multiplicative functions over blocks of $\pi$, namely
$$
\varphi (a_{\pi})= \psi_{\pi_{1}}(1)\psi_{\pi_{2}}(1)\ldots \psi_{\pi_{k}}(1).
$$
Note that in the combinatorics of block depths
given in [AB], the vacuum expectations $\varphi(b_{\pi})$
on the algebra generated by creation and annihilation operators on the interacting Fock space
do not have the above property, i.e.
$$
\varphi(b_{\pi})\neq \beta_{d(1)}\beta_{d(2)}\ldots \beta_{d(k)}
$$
with $\beta_{j}$'s as in (1.5) and $d(j)$'s denoting the depths of $\pi_{j}$'s,
although the moments of the arcsine law are equal to the sums of such products.

Finally, let us remark that Lemmas 10.1-10.2 and Theorem 10.3 can be generalized
to arbitrary $m\in  {\mathbb N}$ by setting ${\rm Inn}(\pi_{k})=0$ for blocks of depth
$d_{k}<m$ and keeping  ${\rm Inn}(\pi_{j})$ as in (5.2) for blocks of depth $d_{j}\geq m$.
This means that up do depth $m$ we have the combinatorics of free probability and
starting fropm depth $m$ we have the combinatorics of monotone probability.
Proofs are very similar but slightly more technical and are omitted.\\
\newpage
\begin{center}
{\sc Bibliography}
\end{center}
[AB] L.~Accardi, M.~Bo\.{z}ejko, Interacting Fock spaces and Gaussianization of
probability measures, {\it Inf. Dim. Anal. Quant. Prob. Rel. Topics} {\bf 1} (1998), 663-670.\\[3pt]
[Av] D.~Avitzour, Free products of $C^{*}$- algebras,
{\it Trans.~Amer.~Math.~Soc.} {\bf 271} (1982), 423-465.\\[3pt]
[BLS] M.Bożejko, M.Leinert, R.Speicher Convolution and limit theorems for conditionally free random variables,
{\it Pacific J.Math.} {\bf 175} (1996), 357-388. \\[3pt]
[BW] M.Bożejko, J.Wysoczański, Remarks on t-transformations of measures and convolutions,
{\it Ann. I.H.Poincare} {\bf 37} (2001), 737-761.\\[3pt]
[FL] U.Franz, R.Lenczewski, Limit theorems for the hierarchy of freeness,
{\it Probab. Math. Stat.} {\bf 19} (1999), 23-41.\\[3pt]
[L1] R.Lenczewski, Unification of independence in quantum probability,
{\it Inf. Dim. Anal. Quant. Prob. Rel. Topics}, {\bf 1} (1998), 383-405.\\[3pt]
[L2] R.Lenczewski, Reduction of free independence to tensor independence,
{\it Inf. Dim. Anal. Quant. Prob. Rel. Topics.} {\bf 7} (2004), 337-360.\\[3pt]
[Lu] Y.~G.~Lu, On the interacting Fock space and the deformed
Wigner law'', {\it Nagoya Math. J.} {\bf 145} (1997), 1-28.\\[3pt]
[M1] N.Muraki, Noncommutative Brownian motion in monotone Fock space, {\it Commun. Math. Phys.}
{\bf 140} (1997), 557-570.\\[3pt]
[M2] N.Muraki, Monotonic independence, monotonic central limit theorem and monotonic law of small numbers,
{\it Inf. Dim. Anal. Quant. Prob. Rel. Topics.} {\bf 4} (2001), 39-58.\\[3pt]
[S] R.~Speicher, A new example of ``independence'' and ``white noise'',
{\it Probab.~Th. Rel.~Fields} {\bf 84} (1990) , 141-159.\\[3pt]
[V] D.Voiculescu, Symmetries of some reduced free product $\mathcal{C}^*$-algebras,
Operator Algebras and their Connections with Topology and Ergodic Theory,
{\it Lecture Notes in Math.} 1132, Springer, Berlin, 1985, 556-588.\\[3pt]
[VDN] D.Voiculescu, K.Dykema, A.Nica, {\it Free Random Variables}, AMS, 1992.
\end{document}